\documentclass[a4paper,11pt]{amsart}
\usepackage{amsmath,amsthm,amssymb}
\usepackage{amscd}
\usepackage{mathrsfs}
\usepackage{graphicx}
\usepackage{xcolor}
\usepackage{url}
\usepackage{array}
\usepackage[all]{xy}
\usepackage{tikz}
\usetikzlibrary{calc}

\usepackage[margin=20truemm]{geometry}

\usepackage{hyperref}

\theoremstyle{definition}
\newtheorem{dfn}{Definition}[section]
\newtheorem{example}[dfn]{Example}
\newtheorem{problem}[dfn]{Problem}
\newtheorem{rem}[dfn]{Remark}

\theoremstyle{plain}

\newtheorem{thm}[dfn]{Theorem}

\newcommand{\Z}{\mathbb{Z}}
\newcommand{\Q}{\mathbb{Q}}
\newcommand{\R}{\mathbb{R}}
\newcommand{\C}{\mathcal{C}}
\newcommand{\gr}{{\rm gr}}
\renewcommand{\L}{\mathcal{L}}
\newcommand{\F}{\mathcal{F}} % N-series filtration
\newcommand{\M}{\mathcal{M}} % mapping class group 
\newcommand{\I}{\mathcal{I}} % Torelli group 
\newcommand{\K}{\mathcal{K}} % Johnson subgroup 
\newcommand{\LC}{{\rm LC}} % lower central series
\newcommand{\J}{{\rm AJ}} % Andreadakis-Johnson filtration
\newcommand{\Aut}{{\rm Aut}}
\newcommand{\Der}{{\rm Der}}
\newcommand{\Hom}{{\rm Hom}}
\newcommand{\Ker}{{\rm Ker}}
\newcommand{\Coker}{{\rm Coker}}
\newcommand{\cok}{\mathfrak{cok}}
\newcommand{\IA}{{\rm IA}}

\newcommand{\ES}{{\rm ES}}
\newcommand{\MC}{{\rm P}\Sigma}
\newcommand{\A}{\mathcal{A}}
\renewcommand{\H}{\mathcal{H}}
\newcommand{\h}{\mathfrak{h}}
\newcommand{\Sp}{{\rm Sp}}
\newcommand{\GL}{{\rm GL}}
\newcommand{\Tr}{{\rm Tr}}
\newcommand{\HS}{\mathcal{HS}}

\tikzset{every picture/.style={line width=1pt}}
\usetikzlibrary{decorations.markings}
\tikzset{->-/.style 2 args={
    postaction={decorate},
    decoration={markings, mark=at position #1 with {\arrow[thick, #2]{>}}}
    },
    ->-/.default={0.5}{}
}
\begin{document}

\title{Around the Andreadakis--Johnson filtration}
\author{Yusuke Kuno}
\address{Department of Mathematics, Tsuda University, 2-1-1 Tsuda-machi, Kodaira-shi, Tokyo 187-8577, Japan}
\email{kunotti@tsuda.ac.jp}
%\date{\today}


\begin{abstract}
The Andreadakis--Johnson filtration and its associated construction, known as the Johnson homomorphism, are useful tools in group theory.
They provide a step-by-step approach to studying the automorphisms of a given group.
After explaining the basics we survey both classical and recent results around the Andreadakis--Johnson filtration, with emphasis on the mapping class group of a once-bordered surface.
\end{abstract}

\maketitle

\section{Introduction}

Commutator calculus is a basic tool in group theory, by which one can capture the non-commutative nature of groups in terms of (graded) Lie algebras.
The theme of this chapter is an application of commutator calculus to the study of automorphisms of groups, which originated with the two mathematicians in the title. 
Let us begin with a brief, non-chronological history.

Let $\Sigma$ be a compact oriented surface of genus $g \ge 1$ with one boundary component, and let $\M$ be the mapping class group of $\Sigma$ relative to the boundary, namely the group of isotopy classes of self-homeomorphisms of $\Sigma$ that fix its boundary pointwise.
For the sake of simplicity, we will confine ourselves to this case.
However, what we present in this chapter admits analogous versions for the mapping class group of other types of surfaces, such as once-punctured surfaces and closed surfaces.

In 1980s, Dennis Johnson wrote a series of papers~\cite{J_aq, J_surv, J_Tor1, J_Tor2, J_Tor3} on the algebraic study of surface automorphisms.
The central topic of Johnson's seminal work was the structure of the Torelli group $\I$, which is an important subgroup of $\M$ defined as the kernel of the action on the first homology of $\Sigma$.
A key ingredient of Johnson's approach is the action of $\M$ on the fundamental group $\pi_1(\Sigma)$, which is free of rank $2g$. 
Using the induced action on the tower of nilpotent quotients of $\pi_1(\Sigma)$ he defined a descending filtration
\[
\M = \M(0) \supset \I = \M(1) \supset \M(2) \supset \cdots  
\]
of $\M$ and introduced an algebraic machinery to capture the associated graded quotient of this filtration.
Johnson showed that for each $k \ge 1$ there is an embedding of the graded piece $\M(k)/\M(k+1)$ into a free abelian group of finite rank, which is nowadays called the $k$th Johnson homomorphism.
Morita~\cite{Mor91, Mor93} refined Johnson's construction, formulated the totality of the Johnson homomorphisms in terms of graded Lie algebras, and posed a question about the characterization of the image. 

Back in 1965, Stylianos Andreadakis wrote an important paper \cite{Andre} on the automorphism group of free groups. 
Let $F_n$ be a free group of rank $n$, and let $\Aut(F_n)$ be the automorphism group of $F_n$.
Andreadakis introduced a descending filtration 
\[
\Aut(F_n) = \A(0) \supset \IA_n = \A(1) \supset \A(2) \supset \cdots
\]
of $\Aut(F_n)$, where the first term $\IA_n$ is called the IA-automorphism group of $F_n$ and defined as the kernel of the action on the abelianization of $F_n$.
In fact, Andreadakis worked more generally with an arbitrary group $G$ and introduced a descending filtration
\[
\Aut(G) = \A_G(0) \supset \A_G(1) \supset \A_G(2) \supset \cdots
\]
of the automorphism group of $G$.  
The mapping class group $\M$ is regarded as a subgroup of $\Aut(F_{2g})$, and the filtration $\{ \A_G(k) \}_{k\ge 0}$ restricts to the filtration that Johnson gave for $\M$.
The filtration $\{ \A_G(k) \}_{k \ge 0}$ is called the Andreadakis-Johnson filtration \cite{CHP11,Satoh16}.
The method used in the definition of the Johnson homomorphism for $\M$ can be applied naturally to $\Aut(G)$.
The resulting object, called the Johnson homomorphism for $\Aut(G)$, serves as a linear approximation to the IA-automorphism group $\IA_G = \A_G(1)$.
 
This chapter is aimed at giving a gentle introduction to the Andreadakis--Johnson filtration and the Johnson homomorphism, with emphasis on the case of the mapping class group $\M$.
Our exposition is not meant to be exhaustive, but we have tried to include some of the most recent results.
There are already good survey articles on this topic, such as \cite{J_surv, Mor99, Satoh16, Hain_surv}.
See also \cite{KK16}.
For more details, we refer the reader to these surveys.

\subsection*{Acknowledgements}
The author would like to thank Athanase Papadopoulos for careful reading of an earlier version of this chapter.
He also would like to express his gratitude to Benson Farb, Richard Hain, Mai Katada, Nariya Kawazumi, Jonathan Pakianathan, Takuya Sakasai, Masatoshi Sato, and Takao Satoh for helpful comments, suggestions and information on relevant references.
This research is supported by JSPS KAKENHI Grant Numbers 23K03121, 24K00520, and 26H01994.

%%%
\section{Andreadakis--Johnson filtration} \label{sec:AJ}

In this section, we review basic materials on the Andreadakis--Johnson filtration.

For a group $G$, we use the following notation.
\begin{itemize}
\item The automorphism group of $G$ is denoted $\Aut(G)$.
\item For $x,y\in G$, their commutator\index{commutator} is $[x,y]:= xyx^{-1}y^{-1}$.
\end{itemize}

Let $R$ be a commutative ring with unit. 
Recall that a Lie algebra\index{Lie algebra} over $R$ is an $R$-module $\mathfrak{g}$ equipped with an $R$-bilinear map $\mathfrak{g} \times \mathfrak{g} \to \mathfrak{g}, (a,b) \mapsto [a,b]$ satisfying the identities $[a,a] = 0$ for any $a \in \mathfrak{g}$ and $[a,[b,c]] + [b, [c,a]] + [c, [a,b]] = 0$ for any $a,b,c \in \mathfrak{g}$.
The latter identity is called the Jacobi identity.
A Lie algebra $\mathfrak{g}$ is called (positively) graded if there is a direct sum decomposition $\mathfrak{g} = \bigoplus_{k \ge 1} \mathfrak{g}_k$ such that $[\mathfrak{g}_k, \mathfrak{g}_l] \subset \mathfrak{g}_{k+l}$ for any $k,l \ge 1$.

Let $\mathfrak{g}$ be a Lie algebra over a commutative ring $R$.
\begin{itemize}
\item A derivation\index{derivation} on $\mathfrak{g}$ is an $R$-linear map $D\colon \mathfrak{g} \to \mathfrak{g}$ such that $D([a,b]) = [D(a),b] + [a,D(b)]$ for any $a,b \in \mathfrak{g}$.
When $\mathfrak{g}$ is graded, $D$ is called of degree $k$ if $D(\mathfrak{g}_l) \subset \mathfrak{g}_{k+l}$ for any $l \ge 1$.

\item The $R$-module $\Der(\mathfrak{g})$ of derivations on $\mathfrak{g}$ forms a Lie algebra over $R$ whose Lie bracket is given by $[D, D'] = D \circ D' - D' \circ D$.
When $\mathfrak{g}$ is positively graded, we consider the Lie algebra $\Der^+(\mathfrak{g}) := \bigoplus_{k \ge 1} \Der^k(\mathfrak{g})$ of positive degree derivations on $\mathfrak{g}$, where $\Der^k(\mathfrak{g})$ is the $R$-module of derivations of degree $k$. 
\end{itemize}

\subsection{The associated Lie algebra of a group} \label{subsec:grG}

We refer to \cite[Chapter 5]{MKS} for more details of the content of this section.

Let $G$ be a group.
The abelianization $G^{\rm abel}$ can be thought of as a first approximation to $G$.
It is an abelian group defined as the quotient of $G$ by the commutator subgroup: $G^{\rm abel}=G/[G,G]$.
Using commutators of higher order, e.g., elements of the form $[x,[y,z]]$, one can obtain a better approximation which extracts the non-commutative nature of $G$.
To formulate it we need to recall the lower central series\index{lower central series} $\{ \Gamma_k G \}_{k \ge 1}$, which is a descending filtration\index{filtration!lower central series} of $G$ by normal subgroups
 \[
 G = \Gamma_1 G \supset \Gamma_2 G \supset \Gamma_3 G \supset \cdots,
 \]
defined inductively by $\Gamma_1 G = G$ and $\Gamma_{k+1} G = [G,\Gamma_k G]$, so that $\Gamma_2 G = [G, G]$, $\Gamma_3 G = [G,[G,G]]$, and so on.
For any $k,l \ge 1$ one has 
$
 [\Gamma_k G, \Gamma_l G] \subset \Gamma_{k+l} G
$.
This implies that for each $k$ the successive quotient
\[
\gr_\LC G(k):= \Gamma_k G/\Gamma_{k+1} G
\]
is abelian.
For $x \in \Gamma_k G$, its coset class in $\gr_\LC G(k)$ is denoted $\{ x\}_k$.
Consider the graded $\Z$-module
\[
\gr_\LC G := \bigoplus_{k \ge 1} \gr_\LC G(k).
\]
By the Hall--Witt identities about group commutators, $\gr_\LC G$ has the structure of a graded Lie algebra over $\Z$ whose Lie bracket is induced from the group commutator in $G$.
Namely, the Lie bracket of the coset classes $\{ x\}_k$ and $\{ y\}_l$ is the coset class of the commutator of $x$ and $y$: 
\[
[ \{ x\}_k, \{ y\}_l ] := \{ [x,y] \}_{k+l}.
\]
For example, the Jacobi identity on $\gr_\LC G$ is a consequence of the following identity: 
\[
[\, {}^z x, [y,z]]\, [\, {}^y z, [x,y]]\, [\, {}^x y, [z,x]] = 1.
\]
Here, we denote ${}^z x = z x z^{-1}$.
The Lie algebra $\gr_\LC G$ is generated by the degree one part $\gr_\LC G(1) =  G^{\rm abel}$, i.e., the iterated Lie bracket
\[
(\gr_\LC G(1))^{\otimes k} \to \gr_\LC G(k), 
\quad 
u_1 \otimes u_2 \otimes \cdots \otimes u_{k-1} \otimes u_k 
\mapsto 
[[\cdots [[u_1, u_2], u_3],\cdots, u_{k-1}], u_k]
\]
is surjective for any $k$. 

When we work over the rationals, we consider $\gr_\LC^\Q G := \gr_\LC G \otimes_\Z \Q$.

\begin{example}
Let $F_n$ be a free\index{free group} group\index{group!free} of rank $n$ with basis $\gamma_1, \ldots, \gamma_n$.
The abelianization $H:=F_n^{\rm abel} = \gr_\LC F_n(1)$ is free abelian of rank $n$, and the coset classes $x_i = \{ \gamma_i \}_i$ constitute a $\Z$-basis for $H$.
The associated graded Lie algebra is the free\index{free Lie algebra} Lie algebra\index{Lie algebra!free} $\L = \L(H) = \bigoplus_{k \ge 1} \L_k$ generated by $H$:
\[
\gr_\LC F_n = \L.
\]
Each graded piece $\gr_\LC F_n(k) = \L_k$ is free abelian of finite rank.
In degrees one and two, we have $\L_1 = H$ and $\L_2 \cong \wedge^2 H$.

There is a nice diagrammatic description of free Lie algebras which holds over a field of characteristic zero.
For definiteness we work over the rationals, and so let $H^\Q = H \otimes_\Z \Q$ and $\L^\Q := \L \otimes_{\Z} \Q$.
A rooted planar trivalent tree is called $H^{\Q}$-labeled if its leaves are labeled by elements of $H^{\Q}$.
Such a tree determines an element of $\L^\Q$ by replacing trivalent vertices with Lie brackets, as in the following example: 
\[
\begin{tikzpicture}[baseline=10pt, x=4mm, y=4mm]
\draw (0,0) node[below=-2pt]{root};
\fill (0,1) circle[radius=1.6pt];
\fill (1,2) circle[radius=1.6pt]; 
\draw (0,1) -- (-2,3) node[above]{$x_1$};
\draw (1,2) -- (0,3) node[above]{$x_2$};
\draw (0,0) -- (0,1) -- (2,3) node[above]{$x_3$};
\end{tikzpicture}
\quad \longmapsto \quad
[x_1, [x_2, x_3]].
\]
This induces a graded $\Q$-linear isomorphism
\[
\L^\Q \overset{\cong}{\longleftarrow} \frac{\Q \{ \text{$H^{\Q}$-labeled rooted planar trivalent trees} \}}{\text{(ML), (AS), (IHX)}}.
\]
Here, (ML) is the multi-linearity relation for the $H^{\Q}$-labelings, (AS) is the anti-symmetry relation about trivalent vertices, namely flipping two edges at a trivalent vertex gives a minus sign, and (IHX)\index{IHX relation} is the following relation:  
\[
\begin{tikzpicture}[baseline=-0.2em, x=1.8em, y=1.8em]
\draw (0,-1)-- (0,1);
\draw (0,1)--(1,1);
\draw (0,1)--(-1,1);
\draw (0,-1)--(1,-1);
\draw (0,-1)--(-1,-1);
\fill (0,-1) circle[radius=1.6pt];
\fill (0,1) circle[radius=1.6pt];
\end{tikzpicture}
\quad - \quad 
\begin{tikzpicture}[baseline=-0.2em, x=1.8em, y=1.8em]
\draw (-1,0)-- (1,0);
\draw (1,0)--(1,1);
\draw (-1,0)--(-1,1);
\draw (1,0)--(1,-1);
\draw (-1,0)--(-1,-1);
\fill (1,0) circle[radius=1.6pt];
\fill (-1,0) circle[radius=1.6pt];
\end{tikzpicture}
\quad + \quad 
\begin{tikzpicture}[baseline=-0.2em, x=1.8em, y=1.8em]
\draw (-1,-1) -- (1,1);
\draw (1,-1) -- (0.2,-0.2);
\draw (-0.2,0.2) -- (-1,1);
\draw (-0.5,-0.5) -- (0.5,-0.5);
\fill (1/2,-0.5) circle[radius=1.6pt];
\fill (-1/2,-0.5) circle[radius=1.6pt];
\end{tikzpicture}
\quad = 0.
\]
\end{example}

\begin{rem}
The above isomorphism for free Lie algebras does not hold over the integers, since there are issues about 2-torsions.
See~\cite{Levine01Add} for more details. 
\end{rem}

There are many other groups whose associated Lie algebra are explicitly known.
We mention only three important examples.
\begin{enumerate}
\item Let $G$ be the fundamental group of a closed orientable surface of genus $g$.
Then $\gr_{\LC} G$ is isomorphic to $\L(H)/\langle \omega \rangle$, where $H$ is the first homology group of the surface and $\langle \omega \rangle$ is the Lie ideal generated by a degree-two element $\omega$ called the symplectic form (see Section~\ref{subsec:review_mcg} for the definition of $\omega$).
This follows from a result of Labute~\cite{Lab70} on one-relator groups. 

\item Let $G$ be the pure braid\index{pure braid group} group\index{group!pure braid} on $n$ strands.
Then $\gr_{\LC} G$ has the following presentation: the generators are degree-one elements $t_{ij} = t_{ji}$ for $1 \le i \neq j \le n$, and the relations are $[t_{ij}, t_{kl}] = 0$ for all distinct indices $i,j,k,l$, and $[t_{ij} + t_{jk}, t_{ik}] = 0$ for all distinct indices $i,j,k$.
This Lie algebra is called the Drinfeld--Kohno Lie algebra \cite{Dri91, Koh85}.

\item Hain~\cite{Hain97, Hain15} gave a presentation for $\gr_{\LC}^{\Q} \I$, where $\I$ is the Torelli group of genus at least $3$.
In fact, Hain gave presentations of the (rational version of the) associated graded Lie algebra of the Torelli group for closed surfaces, as well as for surfaces with punctures and boundaries.
\end{enumerate}

\begin{rem}
An N-series\index{N-series} of a group $G$ is a descending filtration\index{filtration!N-series} $G = \F_1 G \supset \F_2 G \supset \cdots$ such that $[\F_k G, \F_l G] \subset \F_{k+l} G$ for any $k,l \ge 1$.
By Lazard~\cite{Laz54}, the construction of a graded Lie algebra out of a group $G$ applies not only to the lower central series but also to any N-series. 
This generalization is used in the next section to the Andreadakis--Johnson filtration of $\Aut(G)$. 
\end{rem}

\subsection{Andreadakis--Johnson filtration and Johnson homomorphism} \label{subsec:AJJ}

We refer to \cite{Andre}, \cite[Section~3]{Mor91}, \cite[Section~2]{Pa01}, \cite[Section~2]{PaSu12}, \cite{Satoh16} for more details of the content of this section.

Let $G$ be a group.
For $k \ge 1$, the group  $G/\Gamma_{k+1} G$ is called the $k$th nilpotent quotient\index{nilpotent quotient} of $G$.
The action of $\Aut(G)$ on $G$ induces an action on $G/\Gamma_{k+1} G$.

\begin{dfn}
The Andreadakis--Johnson\index{Andreadakis--Johnson filtration} filtration\index{filtration!Andreadakis--Johnson} $\{ \A_G(k) \}_k$ is a descending filtration of $\Aut(G)$ by normal subgroups 
\[
\Aut(G) = \A_G(0) \supset \A_G(1) \supset \A_G(2) \supset \cdots
\]
whose $k$th term is given by 
$\A_G(k):= \Ker \big( \Aut(G) \to \Aut(G/\Gamma_{k+1} G) \big)$.
\end{dfn}

\begin{rem}
If $\bigcap_{k \ge 1} \Gamma_k G$ is trivial (i.e., $G$ is residually nilpotent), $\bigcap_{k \ge 1} \A_G(k)$ is trivial as well.
\end{rem}

By definition, for any $\varphi \in \A_G(k)$ and for any $\gamma \in G$ one has $\varphi(\gamma) \gamma^{-1} \in \Gamma_{k+1} G$.
We denote 
\[
\IA_G := \A_G(1)
\]
and call it the IA-automorphism\index{IA-automorphism group} group\index{group!IA-automorphism} of $G$.
It consists of all automorphisms of $G$ that act on $G^{\rm abel}$ trivially, namely act as the ``Identity Automorphism.''
More generally, for any $\varphi \in \A_G(k)$ and for any $\gamma \in \Gamma_l G$ with $l \ge 1$, it holds that $\varphi(\gamma) \gamma^{-1} \in \Gamma_{k+l} G$. 
This key observation together with a Lie algebra theoretic formulation results in an injective graded Lie algebra homomorphism
\[
\tau \colon \gr_\J \IA_G \to \Der^+(\gr_\LC G).
\]
We call $\tau$ the Johnson\index{Johnson homomorphism} homomorphism\index{homomorphism!Johnson} and its degree $k$ component the $k$th Johnson homomorphism.
Let us explain the definition of $\tau$ in more detail.

\begin{itemize}

\item
{\em The source.}
The filtration $\{ \A_G(k) \}_{k \ge 1}$ is an N-series, i.e., for any $k,l \ge 1$ one has 
$
[\A_G(k), \A_G(l)] \subset \A_G(k+l)
$.
Thus the successive quotient
\[
\gr_\J \IA_G(k) := \A_G(k)/\A_G(k+1)
\]
is abelian and the graded $\Z$-module
\[
\gr_\J \IA_G:= \bigoplus_{k \ge 1} \gr_\J \IA_G(k)
\]
has the structure of a graded Lie algebra over $\Z$, whose Lie bracket is induced from the group commutator in $\IA_G$.

\item
{\em The target.}
Let $\Der^k(\gr_\LC G)$ be the $\Z$-module of derivations on $\gr_\LC G$ which are of degree $k$.
The direct sum $\Der^+(\gr_\LC G) = \bigoplus_{k \ge 1} \Der^k(\gr_\LC G)$ is the Lie algebra of positive degree derivations on $\gr_\LC G$.

\item
{\em The map $\tau$.}
For any $\varphi \in \A_G(k)$, the map
\[
\tau_k(\varphi) \colon \gr_\LC G \to \gr_\LC G, 
\quad \{ \gamma \}_l \mapsto \{ \varphi(\gamma) \gamma^{-1} \}_{k+l}
\]
is a well-defined graded $\Z$-module map of degree $k$. 
Furthermore, $\tau_k(\varphi)$ is a derivation on $\gr_\LC G$.
Thus one gets a well-defined map
\[
\tau_k \colon \A_G(k) \to \Der^k(\gr_\LC G), 
\quad \varphi \mapsto \tau_k(\varphi).
\]
This is a group homomorphism, and by construction, its kernel is $\A_G(k+1)$. 
We abuse the same letter for the induced injective homomorphism 
$
\tau_k \colon \gr_\J \IA_G(k) \to 
\Der^k(\gr_\LC G)
$.
The collection $\tau = (\tau_k)_{k \ge 1}$ defines an injective graded $\Z$-module map from $\gr_\J \IA_G$ to $\Der^+(\gr_{\LC} G)$.
The claim that $\tau$ is a Lie algebra homomorphism is equivalent to the following: 
for any $\varphi \in \A_G(k)$ and $\psi \in \A_G(l)$, it holds that $[\varphi, \psi] \in \A_G(k+l)$ and 
\[
\tau_{k+l} ([\varphi, \psi]) 
= \tau_k(\varphi) \circ \tau_l(\psi) - \tau_l(\psi) \circ \tau_k(\varphi).
\]
\end{itemize}

The group $\Aut(G)/\IA_G$ is naturally isomorphic to the group of automorphisms of $G^{\rm abel}$ induced from automorphisms of $G$. 
This group acts equivariantly on the source and target of $\tau_k$.
The actions are given as follows:
\begin{itemize}
\item 
{\em On the source.}
The group $\Aut(G) = \A_G(0)$ acts on $\A_G(k)$ by conjugation and hence on the quotient $\gr_\J \IA_G (k) = \A_G(k)/\A_G(k+1)$ as well.
The action of $\IA_G = \A_G(1)$ is trivial.
Thus the quotient $\Aut(G)/\IA_G$ acts on $\gr_\J \IA_G(k)$. 
\item
{\em On the target.}
The action of $\Aut(G)$ on the lower central series of $G$ induces an action on $\gr_\LC G$.
Therefore $\Aut(G)$ acts on the $\Z$-module $\Der^k(\gr_\LC G)$ by conjugation: for $D \in \Der^k(\gr_\LC G)$, the derivation $\varphi \cdot D$ is defined by $(\varphi \cdot D)(u) = \varphi(D (\varphi^{-1}(u)))$.
\end{itemize}

When we work over the rationals, we consider 
\[
\tau^\Q
\colon \gr_{\J}^{\Q} \IA_G \to \Der^+( \gr_{\LC}^{\Q} G), 
\]
where $\gr_{\J}^{\Q} \IA_G = \gr_{\J} \IA_G \otimes_{\Z} \Q$ and $\tau^\Q = \tau \otimes_{\Z} {\rm id}_{\Q}$.
Passing to the rational setting not only eliminates difficulties with torsion but also makes representation theoretic considerations easier. 

As the Johnson homomorphism is injective by construction, the natural problem is the following.

\begin{problem}[The Johnson cokernel problem] \label{pb:JCP}
Given a group $G$, determine the cokernel of the Johnson homomorphism $\tau$ as well as its rational version $\tau^\Q$.
\end{problem}

\subsection{The automorphisms of free groups} \label{subsec:autFn}

Let us consider the case where $G$ is a free group, $G = F_n$.
To simplify notation we write  $\IA_{F_n} = \IA_n$.
It is well known that any automorphism of $H = F_n^{\rm abel}$ is induced from automorphisms of $F_n$, namely we have $\Aut(F_n)/\IA_{n} \cong \Aut(H)$.
Thus we have the following exact sequence:
\[
1 \to \IA_n \to \Aut F_n \to \Aut(H) \to 1.
\]
The Johnson homomorphism for $\IA_n$ is of the form 
\[
\tau \colon \gr_\J \IA_{n} \to \Der^+(\L), 
\]
which is $\Aut(H)$-equivariant.
In the rational setting, we have 
\[
\tau^\Q \colon \gr_{\J}^{\Q} \IA_n \to \Der^+(\L^\Q),
\]
and it is known that this becomes a map of $\Aut(H^{\Q})$-modules. 
See \cite{Patzt} for instance.
Note that the choice of a basis for $H$ defines isomorphisms $\Aut(H) \cong \GL_n(\Z)$ and $\Aut(H^{\Q}) \cong \GL_n(\Q)$.

For the automorphism group of free groups, Problem~\ref{pb:JCP} has been settled in the stable range.
To begin with, we explain the case $k=1$.
The first Johnson homomorphism is of the form
\[
\tau_1\colon \IA_n \to \Der^1(\L) = \Hom(H,\L_2), 
\]
and it is an $\Aut(H)$-equivariant surjection.
This fact was independently proved by Formanek~\cite{For90},  Cohen and Pakianathan~\cite{CohPak}, Farb~\cite{Farb} and Kawazumi~\cite{Kaw05}.
It also holds that $\tau_1$ induces an $\Aut(H)$-equivariant isomorphism $\IA_n^{\rm abel} \overset{\cong}{\to} \Hom(H,\L_2)$ and $\Gamma_2 \IA_n = \A_{F_n}(2)$.

\begin{rem}
The description of $\IA_n^{\rm abel}$ as an abelian group was first given by Bachmuth~\cite[Theorem on p.\ 7]{Bach}.
 Andreadakis~\cite[Sections~4 and 5]{Andre} also gave it implicitly and independently.
\end{rem}

To explain the case $k\ge 2$, let us recall the classical fact (e.g., see \cite[Section~5.6]{MKS}) that the degree $k$ part of the free Lie algebra $\L_k$ is naturally embedded into the $k$-fold tensor product $H^{\otimes k}$.
This is done by reinterpreting the Lie bracket symbol as the algebra commutator in the tensor algebra generated by $H$: for example, $\L_2 \ni [x_i, x_j] \mapsto x_i \otimes x_j - x_j \otimes x_i \in H^{\otimes 2}$.
Let $\C_n(k)$ be the quotient of $H^{\otimes k}$ by the cyclic permutations of the components. 
Elements of $\C_n(k)$ can be seen as linear combinations of cyclic words in $\{ x_i \}_i$ of length $k$.
It is sometimes convenient to use the following notation for cyclic words: a diagrammatic notation and an absolute value notation.
\[
\begin{tikzpicture}[baseline=0pt, x=4mm, y=4mm]
\draw (0,0) circle[radius=2];
\draw (0:2) node[right]{$x_1$};
\fill (0:2) circle[radius=1.6pt];
\draw (60:2) node[above right=-2pt]{$x_1$};
\fill (60:2) circle[radius=1.6pt];
\draw (120:2) node[above left=-2pt]{$x_3$};
\fill (120:2) circle[radius=1.6pt];
\draw (180:2) node[left]{$x_2$};
\fill (180:2) circle[radius=1.6pt];
\draw (240:2) node[below left=-2pt]{$x_1$};
\fill (240:2) circle[radius=1.6pt];
\draw (300:2) node[below right=-2pt]{$x_4$};
\fill (300:2) circle[radius=1.6pt];
\end{tikzpicture}
\quad 
\longleftrightarrow
\quad
|x_1 x_1 x_3 x_2 x_1 x_4|
= |x_1 x_3 x_2 x_1 x_4 x_1|.
\]

Satoh~\cite{Satoh06, Satoh12} introduced the following map:
\[
\Phi_k \colon \Der^k(\L) \overset{\cong}{\to} 
H^* \otimes \L_{k+1} \to
H^* \otimes H^{\otimes (k+1)}
\to
H^{\otimes k} 
\to
\C_n(k).
\]
Here, the first map is the restriction $D \mapsto D|_H \in \Hom(H, \L_{k+1}) = H^* \otimes \L_{k+1}$, the second map is induced by the embedding $\L_{k+1} \hookrightarrow H^{\otimes (k+1)}$, the third map is the contraction of the first and the second components (namely it maps $\varphi \otimes x_1 \otimes x_2 \otimes \cdots \otimes x_{k+1}$ to $\varphi(x_1)\, x_2 \otimes \cdots \otimes x_{k+1}$), and the last map is the canonical projection.

\begin{thm}[Satoh, Darn\'{e}] \label{thm:JCIAn}
For $k \ge 2$ and $n \ge k+2$, the following sequence is exact:
\[
0 \to \gr_\J \IA_n(k) \overset{\tau_k}{\to}
\Der^k(\L) \overset{\Phi_k}{\to}
\C_n(k) \to 0,
\] 
therefore $\Coker (\tau_k) \cong \C_n(k)$ in this range.
\end{thm}

This result follows directly by combining results of Satoh~\cite[Theorem~3.1]{Satoh12} and Darn\'{e}~\cite{Dar19} (see Theorem~\ref{thm:stbAP} below for Darn\'{e}'s result).
See also \cite[Section~2.5.3]{Dar19}.
Over the rationals, Massuyeau and Sakasai~\cite{MasSak} showed this exact sequence in the same stable range, based on the above-mentioned result of Satoh.
The fact that $\Phi_k$ vanishes on the image of $\tau_k$ can be proved by several ways.
For one proof, see Massuyeau and Sakasai~\cite[Proposition~5.3]{MasSak}.

%%%
\subsection{Restriction to subgroups}

If $\H$ is a subgroup of $\Aut(G)$, then we can restrict the Andreadakis--Johnson filtration to $\H$ by setting
$\H(k):= \H \cap \A_G(k)$.
This gives an N-series on $\H$:
\[
\H = \H(0) \supset \H(1) \supset \H(2) \supset \cdots.
\]
We call $\I_\H = \H(1) = \H \cap \IA_G$ the Torelli part of $\H$.
One can form the associated graded Lie algebra $\gr_\J \I_\H:= \bigoplus_{k \ge 1} \H(k)/\H(k+1)$, and the Johnson homomorphism for $\H$ is defined by restriction: 
it is an injective graded Lie algebra homomorphism of the form
\[
\tau \colon \gr_\J \I_\H \to \Der^+(\gr_\LC G).
\]
In this case, $\tau$ is equivariant under the action of $\H/\I_\H$. 
We also consider the rational version
\[
\tau^\Q \colon \gr^{\Q}_\J \I_{\H} \to \Der^+(\gr^{\Q}_{\LC} G).
\]

\begin{rem}
One may find a naturally defined Lie subalgebra of the target of $\tau$ arising from some characteristic property of $\H$.
In this case, it is natural to think of this Lie subalgebra as a new target of $\tau$ when considering the Johnson cokernel problem. 
Indeed, this is the case for the mapping class group as we will explain in Sections~\ref{subsec:MR} and \ref{subsec:JCP}.
\end{rem}

%%%
\subsection{The Andreadakis problem} \label{subsec:Ap}

The Andreadakis--Johnson filtration and the Johnson homomorphism give step-by-step approximations to the IA-automorphism group $\IA_G$ using its action on $G$. 
However, one could work more directly using the group structure of $\IA_G$ itself, namely the lower central series of $\IA_G$ and the associated graded Lie algebra $\gr_\LC \IA_G$.

Comparison of the two approaches has been of interest since the work of Andreadakis.
In fact, they are intimately related.
We have $\Gamma_1 \IA_G = \A_G(1)$ by definition, and by the N-series property of the Andreadakis--Johnson filtration we have $\Gamma_k \IA_G \subset \A_G(k)$ for any $k \ge 1$.
If $\H$ is a subgroup of $\Aut(G)$, then $\Gamma_k \I_\H \subset \H(k)$ for any $k \ge 1$.

\begin{problem} \label{pb:Andre}
Is $\Gamma_k \IA_G = \A_G(k)$?\index{Andreadakis problem}
Given a subgroup $\H$ of  $\Aut(G)$, is $\Gamma_k \I_\H = \H(k)$?
\end{problem}

It is Andreadakis who first studied this problem extensively when $G = F_n$. 
In \cite[Section~6]{Andre}, he proved that $\Gamma_k \IA_2 = \A_{F_2}(k)$ for $n=2$ and all $k$, and conjectured that $\Gamma_k \IA_n = \A_{F_n}(k)$ for all $n$ and $k$.
For $k=2$, Bachmuth~\cite[Lemma~5]{Bach} first proved $\Gamma_2 \IA_n = \A_{F_n}(2)$ for all $n \ge 2$. 
As we have explained in Section~\ref{subsec:autFn}, this result also follows from the computation of $\tau_1$.
For $k=3$, Satoh~\cite{Satoh19} proved that $\Gamma_3 \IA_n = \A_{F_n}(3)$ for all $n \ge 3$, improving an earlier result of Pettet~\cite{Pettet}. 

Notably, Bartholdi~\cite{Barth, BarthEr} showed that $\A_{F_3}(4)/\Gamma_4 \IA_3$ is a nontrivial torsion group and $\A_{F_3}(5)/\Gamma_5 \IA_3$ has a nontrivial free part.
Thus the original Andreadakis conjecture is disproved for $n = 3$, and it remains open for $n \ge 4$. 
One can ask Problem~\ref{pb:Andre} for subgroups of $\Aut(F_n)$.
Several interesting results have been obtained.
On the one hand, there are groups for which Problem~\ref{pb:Andre} have affirmative answers: 
\begin{enumerate}
\item The (pure) braid group on $n$-strands, viewed as a subgroup of $\Aut(F_n)$ by the Artin representation. See \cite[Theorem~6.2 and Remark~6.4]{Dar21}.
\item The group of triangular automorphisms of $F_n$. See \cite{Satoh17}, \cite[Theorem~5.4]{Dar21}.
\end{enumerate}
See also Darn\'{e}~\cite{Dar25} for another affirmative example from the mapping class group of punctured spheres.
On the other hand, there is a group for which the answer is negative.
The McCool\index{McCool group} group\index{group!McCool} $\MC_n$, also known as the pure loop braid group and the pure welded braid group, is defined as the group of $\varphi \in \Aut(F_n)$ such that $\varphi(\gamma_i)$ is conjugate to $\gamma_i$ for every $i = 1, \ldots, n$.  
Darn\'{e}, Enomoto and Satoh~\cite{DES} showed that the McCool group does not possesses the Andreadakis property by proving that the quotient $\MC_n(7)/\Gamma_7 \MC_n$ contains a free abelian group of rank $\binom{n}{3}$ for $n \ge 3$.

The inclusion $\Gamma_k \IA_G \subset \A_G(k)$ yields the comparison homomorphism
\[
c \colon \gr_\LC \IA_G \to \gr_\J \IA_G
\]
and its rational version $c^{\Q} \colon \gr_{\LC}^{\Q} \IA_G \to \gr_{\J}^{\Q} \IA_G$.
If $\H$ is a subgroup of $\Aut(G)$, we similarly have 
\[
c \colon \gr_\LC \I_\H \to \gr_\J \I_\H
\]
as well as the rational version $c^{\Q} \colon \gr^{\Q}_\LC \I_\H \to \gr^{\Q}_\J \I_\H$.
The associated graded version of the Andreadakis problem is the following.

\begin{problem}
Is $c$ an isomorphism? Is $c^{\Q}$ an isomorphism?
\end{problem}

For the automorphism group of free groups, the following strong result is known.

\begin{thm}[Darn\'{e} \cite{Dar19}] \label{thm:stbAP}
If $n \ge k+2$, then 
$c\colon \gr_{\LC} \IA_n (k) \to \gr_{\J} \IA_n(k)$ is surjective.
\end{thm}
 
\begin{rem}
We can use the graded Lie algebra $\gr_{\LC} \IA_G$ associated with the lower central series of $\IA_G$ as an approximation to $\IA_G$.
Why do we consider the Johnson filtration and the associated graded Lie algebra $\gr_\J \IA_G$? 
The Johnson homomorphism is ``designed to be injective'', and it gives a natural ambient space of $\gr_\J \IA_G$, namely $\Der^+(\gr_{\LC} G)$. 
One could view this as an advantage of the Andreadakis--Johnson filtration over the lower central series of $\IA_G$.
Just the existence of the Johnson homomorphism tell us some information on the Andreadakis--Johnson filtration.
For example, for the automorphism group of a free group, each graded piece $\gr_\J \IA_n(k)$ turns out to have no torsion.
This is because $\Der^k (\L) \cong \Hom(H, \L_{k+1})$ has no torsion. 
The same is true also for any subgroup of $\Aut(F_n)$.
On the other hand, $\gr_{\LC} \IA_n(k)$ may have torsion.
\end{rem}

\subsection{Finite generation} \label{subsec:fg}

One of the basic problems for a given group is to determine whether it is finitely generated.

\begin{problem}
Determine whether $\A_G(k)$ and $\Gamma_k \IA_G$ are finitely generated.
More generally, for a subgroup $\H$ of $\Aut(G)$, determine whether $\H(k)$ and $\Gamma_k \I_\H$ are finitely generated.
\end{problem}

For the automorphism group of a free group, it is well known that both $\A_{F_n}(0) = \Aut(F_n)$ and $\A_{F_n}(1) = \IA_n$ are finitely generated.
See \cite[Section~3.5]{MKS} for explicit generators.
For some time it remained unknown whether $\A_{F_n}(2)$ was finitely generated.
Ershov and He~\cite{EH18} proved that for $n \ge 4$ any subgroup of $\IA_n$ containing $\A_{F_n}(2)$ is finitely generated.
This was generalized to the higher terms:

\begin{thm}[\text{Church--Ershov--Putman \cite[Theorem~E]{CEP22}}]
For $k \ge 3$ and $n \ge 2k+3$, any subgroup of $\IA_n$ containing $\A_{F_n}(k)$ is finitely generated. 
\end{thm}

Church, Ershov and Putman~\cite[Theorem~D]{CEP22} also obtained a similar result on the finite generation of the lower central series of $\IA_n$.

This finiteness property does not always hold for subgroups of $\Aut(F_n)$.
A counter-example is the McCool group $\MC_n$.
Satoh~\cite{Satoh12MC} proved that for $n \ge 3$ and $k \ge 2$ both $\MC_n(k)$ and $\Gamma_k \MC_n$ are not finitely generated.

%%%%%%%%%%%%%
%%%%%%%%%%%%%
\section{The mapping class group} \label{sec:mcg}

In this section, we survey results on the Johnson homomorphism for the mapping class group of a once-bordered surface.
We refer to \cite{FM12} for more details on the basics of the mapping class group. 

\subsection{A quick review of the mapping class group} \label{subsec:review_mcg}
Let $g \ge 1$, and let $\Sigma$ be a compact oriented surface of genus $g$ with one boundary component.
We choose a basepoint from the boundary of $\Sigma$ and let $\pi = \pi_1(\Sigma,*)$ be the fundamental group of $\Sigma$.
The group $\pi$ is a free group of rank $2g$.
It has a distinguished element $\zeta$ that corresponds to the boundary of $\Sigma$.
One can take a standard basis $\{ \alpha_i, \beta_i \}_{i=1}^g$ for $\pi$, as shown in the following figure.
Then we have $\zeta = \prod_{i=1}^g [\alpha_i, \beta_i]$.

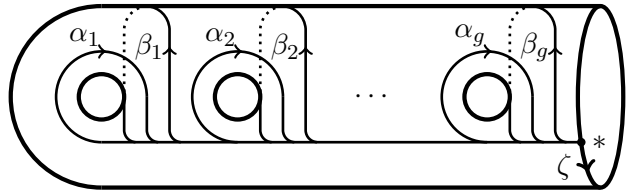
\[
\begin{tikzpicture}[x=3mm, y=3mm]
\draw[ultra thick] (0,0) circle[x radius=1, y radius=4];
\draw[ultra thick] (0,4) -- (-22,4);
\draw[ultra thick] (0,-4) -- (-22,-4);
\draw[ultra thick] (-22,4) arc[x radius=4, y radius=4, start angle=90, end angle=270];
\draw[ultra thick] (-5,0) circle[radius=1];
\draw (-10,0) node{$\cdots$};
\draw[ultra thick] (-16,0) circle[radius=1];
\draw[ultra thick] (-22,0) circle[radius=1];]
\fill (-0.85,-2) circle[radius=1.6pt] node[right]{$*$};
\draw[ultra thick, ->] (-0.6,-3.1) -- (-0.55,-3.3);
\draw (-0.55,-3.1) node[left=2pt]{$\zeta$};
\draw (-22,-2) -- (-0.85,-2);
\draw (-20,0) arc[radius=2, start angle=0, end angle=270];
\draw (-20,0) -- (-20,-1.5) to[bend right=40] (-19.5,-2);
\draw (-18.5,-2) to[bend left=40] (-19,-1.5) -- (-19,3) to[bend right=40] (-20,4);
\draw[dotted] (-20,4) to[bend right=40] (-21,3) -- (-21,0);
\draw (-21,0) -- (-21,-1.5) to[bend right=40] (-20.5,-2);
\draw[->] (-22,2) -- (-21.8,2) node[above left=-2pt]{$\alpha_1$};
\draw[->] (-19,2) -- (-19,2.2) node[left=-2pt]{$\beta_1$};
\draw (-14,0) arc[radius=2, start angle=0, end angle=270];
\draw (-14,0) -- (-14,-1.5) to[bend right=40] (-13.5,-2);
\draw (-12.5,-2) to[bend left=40] (-13,-1.5) -- (-13,3) to[bend right=40] (-14,4);
\draw[dotted] (-14,4) to[bend right=40] (-15,3) -- (-15,0);
\draw (-15,0) -- (-15,-1.5) to[bend right=40] (-14.5,-2);
\draw[->] (-16,2) -- (-15.8,2) node[above left=-2pt]{$\alpha_2$};
\draw[->] (-13,2) -- (-13,2.2) node[left=-2pt]{$\beta_2$};
\draw (-3,0) arc[radius=2, start angle=0, end angle=270];
\draw (-3,0) -- (-3,-1.5) to[bend right=40] (-2.5,-2);
\draw (-1.5,-2) to[bend left=40] (-2,-1.5) -- (-2,3) to[bend right=40] (-3,4);
\draw[dotted] (-3,4) to[bend right=40] (-4,3) -- (-4,0);
\draw (-4,0) -- (-4,-1.5) to[bend right=40] (-3.5,-2);
\draw[->] (-5,2) -- (-4.8,2) node[above left=-2pt]{$\alpha_g$};
\draw[->] (-2,2) -- (-2,2.2) node[left=-2pt]{$\beta_g$};
\end{tikzpicture}
\]

We use a shorthand notation $H = \pi^{\rm abel}$.
The symbol $H$ stands for ``homology'', as we have a canonical isomorphism $H \cong H_1(\Sigma; \Z)$.
The intersection pairing on $\Sigma$ is a skew-symmetric non-degenerate pairing $(\ \cdot \ ) \colon H \times H \to \Z$. 
The coset classes $a_i = \{ \alpha_i \}_1$ and $b_i = \{ \beta_i \}_1$ constitute a symplectic basis\index{symplectic basis} for $(H, (\cdot, \cdot))$, i.e.,  they satisfy $(a_i \cdot b_j) = \delta_{ij}$ and $(a_i \cdot a_j) = (b_i \cdot b_j) = 0$.

The element $\zeta$ belongs to $\Gamma_2 \pi = [\pi, \pi]$, and 
it projects to a degree-two element $\omega = \{ \zeta \}_2$ in the free Lie algebra $\L = \L(H) = \gr_{\LC} \pi$.
Explicitly, we have $\omega = \sum_{i=1}^g [a_i, b_i]$.
We call $\omega$ the symplectic form.\index{symplectic form}

Let $\M$ be the mapping class\index{mapping class group} group\index{group!mapping class} of $\Sigma$, namely the group of self-homeomorphisms of $\Sigma$ modulo isotopies.
Here, the self-homeomorphisms and the isotopies that we consider fix the boundary of $\Sigma$ pointwise. 
Therefore, $\M$ acts naturally on the free group $\pi$.

\begin{thm}[the Dehn--Nielsen theorem]
The action of $\M$ on $\pi$ induces an injective homomorphism
\[
\rho \colon \M \to \Aut(\pi).
\]
The image coincides with automorphisms $\varphi \in \Aut(\pi)$ such that $\varphi(\zeta) = \zeta$. \index{theorem!Dehn--Nielsen}
\end{thm}

Through the Dehn--Nielsen map $\rho$ we can regard $\M$ as a subgroup of $\Aut(\pi)$, and we obtain the Andreadakis--Johnson filtration on $\M$:
\[
\M = \M(0) \supset \M(1) \supset \M(2) \supset \cdots.
\]
Here, $\M(k) = \M \cap \A_{\pi}(k)$.
The group $\I = \I_{\M} = \M(1)$ is called the Torelli\index{Torelli group} group,\index{group!Torelli} and the next term $\K = \M(2)$ the Johnson kernel.\index{Johnson kernel}
The filtration $\{ \M(k) \}_{k}$ is traditionally called the Johnson\index{Johnson filtration} filtration.\index{filtration!Johnson} 

There are beautiful descriptions of the first three terms of the Andreadakis--Johnson filtration on $\M$ by Dehn twists.
Let $\gamma$ be a simple closed curve in $\Sigma$.
Note that a closed tubular neighborhood of $\gamma$ is an annulus, since $\Sigma$ is orientable.
The Dehn twist\index{Dehn twist} $t_{\gamma}$ along $\gamma$ is a self-homeomorphism of $\Sigma$ defined by cutting $\Sigma$ along $\gamma$, rotating it, and gluing it back:
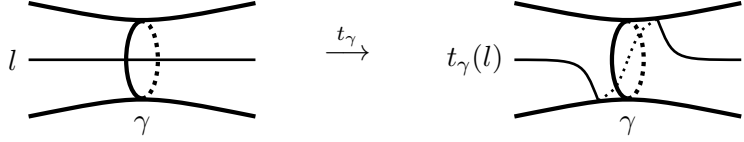
\[
\begin{tikzpicture}[baseline=0pt, x=3mm, y=3mm]
\draw[ultra thick] (-5,2.5) ..controls(0,1.5).. (5,2.5);
\draw[ultra thick] (-5,-2.5) ..controls(0,-1.5).. (5,-2.5);
\draw[ultra thick] (0,1.7) arc[x radius=0.7, y radius=1.7, start angle=90, end angle=270]; 
\draw[ultra thick, dotted] (0,-1.7) arc[x radius=0.7, y radius=1.7, start angle=-90, end angle=90]; 
\draw (-5,0) -- (5,0);
\draw (0,-2) node[below]{$\gamma$};
\draw (-5,0) node[left]{$l$};
\end{tikzpicture}
\qquad 
\overset{t_\gamma}{\longrightarrow}
\qquad 
\begin{tikzpicture}[baseline=0pt, x=3mm, y=3mm]
\draw[ultra thick] (-5,2.5) ..controls(0,1.5).. (5,2.5);
\draw[ultra thick] (-5,-2.5) ..controls(0,-1.5).. (5,-2.5);
\draw[ultra thick] (0,1.7) arc[x radius=0.7, y radius=1.7, start angle=90, end angle=270]; 
\draw[ultra thick, dotted] (0,-1.7) arc[x radius=0.7, y radius=1.7, start angle=-90, end angle=90]; 
\draw (-5,0) ..controls(-2,0).. (-1.3,-1.8);
\draw (5,0) ..controls(2,0).. (1.3,1.8);
\draw[dotted] (-1.3,-1.8) ..controls(-0.5,-1.3).. (0,0);
\draw[dotted] (1.3,1.8) ..controls(0.5,1.3).. (0,0);
\draw (0,-2) node[below]{$\gamma$};
\draw (-5,0) node[left]{$t_{\gamma}(l)$};
\end{tikzpicture}
\]

\begin{enumerate}

\item[(0)]
The mapping class group $\M$ is generated by Dehn twists.
What is more, there are explicit sets of generators of $\M$ consisting of finitely many Dehn twists.

\item[(1)] 
The Torelli group $\I$ is generated by BP-maps\index{BP-map} if $g \ge 3$.
Here, a pair of disjoint simple closed curves $(\gamma, \delta)$ is called a bounding pair if their union is the boundary of a subsurface of $\Sigma$; see below.
The product $t_\gamma t_{\delta}^{-1}$ is called a BP-map.
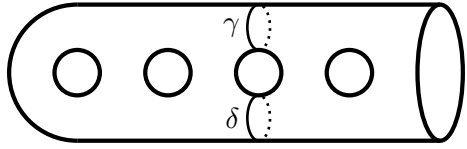
\[
\begin{tikzpicture}[x=3mm, y=3mm]
\draw[ultra thick] (0,0) circle[x radius=1, y radius=3];
\draw[ultra thick] (0,3) -- (-16,3);
\draw[ultra thick] (0,-3) -- (-16,-3);
\draw[ultra thick] (-16,3) arc[x radius=3, y radius=3, start angle=90, end angle=270];
\draw[ultra thick] (-4,0) circle[radius=1];
\draw[ultra thick] (-8,0) circle[radius=1];
\draw[ultra thick] (-12,0) circle[radius=1];
\draw[ultra thick] (-16,0) circle[radius=1];]
\draw (-8,3) arc[x radius=0.5, y radius=1, start angle=90, end angle=270];
\draw[dotted] (-8,1) arc[x radius=0.5, y radius=1, start angle=-90, end angle=90];
\draw (-8,-1) arc[x radius=0.5, y radius=1, start angle=90, end angle=270];
\draw[dotted] (-8,-3) arc[x radius=0.5, y radius=1, start angle=-90, end angle=90];
\draw (-8,2) node[left=3pt]{$\gamma$};
\draw (-8,-2) node[left=3pt]{$\delta$};
\end{tikzpicture}
\]
A remarkable result of Johnson~\cite{J_Tor1} is that the Torelli group is finitely generated.
In fact, he gave a generating set consisting of finitely many BP-maps.

\item[(2)]
Johnson~\cite{J_Tor2} proved that the group $\K$ is generated by BSCC-maps.\index{BSCC-map}
Here, a simple closed curve $\gamma$ is called separating if it bounds a subsurface of $\Sigma$; see below.
The Dehn twist along $\gamma$ is called a BSCC-map (bounding simple closed curves).
See also Putman~\cite{Put18} for another proof and a generalization.
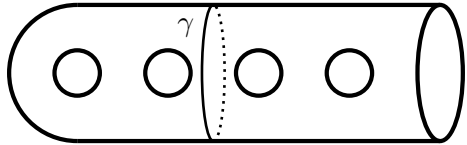
\[
\begin{tikzpicture}[x=3mm, y=3mm]
\draw[ultra thick] (0,0) circle[x radius=1, y radius=3];
\draw[ultra thick] (0,3) -- (-16,3);
\draw[ultra thick] (0,-3) -- (-16,-3);
\draw[ultra thick] (-16,3) arc[x radius=3, y radius=3, start angle=90, end angle=270];
\draw[ultra thick] (-4,0) circle[radius=1];
\draw[ultra thick] (-8,0) circle[radius=1];
\draw[ultra thick] (-12,0) circle[radius=1];
\draw[ultra thick] (-16,0) circle[radius=1];]
\draw (-10,3) arc[x radius=0.5, y radius=3, start angle=90, end angle=270];
\draw[dotted] (-10,-3) arc[x radius=0.5, y radius=3, start angle=-90, end angle=90];
\draw (-10,2) node[left=3pt]{$\gamma$};
\end{tikzpicture}
\]
\end{enumerate}

\begin{rem}
The size of the BP-map generators that Johnson gave for $\I$ has order $O(2^g)$.
Putman~\cite{Put12} gave another explicit finite generating set of $\I$ which consists of BP-maps and BSCC-maps and whose cardinality is much smaller.
\end{rem}

Let $\Sp(H)$ be the group of automorphisms of $H$ that preserve the intersection pairing $(\cdot, \cdot)$.
Any automorphism of $H$ that is induced from $\M$ must preserve $(\cdot, \cdot)$.
In fact, we have $\M/\I \cong \Sp(H)$.
Thus we have the following exact sequence:
\[
1 \to \I \to \M \to \Sp(H) \to 1.
\]
Note that the choice of a symplectic basis for $(H, (\cdot,\cdot))$ defines an isomorphism $\Sp(H) \cong \Sp_{2g}(\Z)$.

\begin{rem}
When $g=1$, the above exact sequence becomes a central extension
\[
0 \to \Z \to \M \to {\rm SL}_2(\Z) \to 1.
\]
Let $\gamma$ and $\delta$ be non-separating simple closed curves in $\Sigma$ which intersect transversely in one point.
Then the Dehn twists $t_{\gamma}$ and $t_{\delta}$ generate $\M$, and we have the presentation $\M = \langle t_{\gamma}, t_{\delta} \mid t_{\gamma} t_{\delta} t_{\gamma} = t_{\delta} t_{\gamma} t_{\delta} \rangle$.
Hence, $\M$ is isomorphic to the braid group on $3$-strands.
The kernel of $\M \to {\rm SL}_2(\Z)$ is generated by $(t_{\gamma}t_{\delta})^6$.
In particular, the lower central series of the Torelli group is given by $\Gamma_1\I = \Z$ and $\Gamma_k \I = \{ e\}$ for $k \ge 2$.
As for the Johnson filtration, we have $\M(1) = \M(2) = \Z$ and $\M(k) = \{ e\}$ for $k \ge 3$.
\end{rem}

%%%
\subsection{Morita's refinement} \label{subsec:MR}

A priori, the Johnson homomorphism for the mapping class group $\M$ takes values in $\Der^+(\L)$, which is the same as the target of the Johnson homomorphism for $\Aut(\pi)$. 
Since $\pi$ is the fundamental group of the once-bordered surface $\Sigma$, the free Lie algebra $\L = \gr_{\LC} \pi$ has a distinguished element $\omega$ which corresponds to the boundary element $\zeta$ in $\pi$.
Morita \cite[Section~3]{Mor93} showed that the boundary condition $\varphi(\zeta) = \zeta$ imposes a natural restriction on the image of $\tau$: it holds that $\tau_k(\varphi)(\omega) = 0$ for any $\varphi \in \M(k)$.
A derivation $D \in \Der^+(\L)$ is called symplectic\index{symplectic derivation} if $D(\omega) = 0$.\index{derivation!symplectic}
Morita introduced the space 
\[
\h:= \{ u \in \Der^+(\L) \mid D(\omega) = 0 \}
\]
of symplectic derivations on $\L$.
This is a graded Lie subalgebra of $\Der^+(\L)$ and closed under the action of $\Sp(H)$. 
Thus, for the mapping class group it is more natural to regard $\tau$ as taking values in $\h$:
\[
\tau \colon \gr_\J \I \to \h.
\]
This is an $\Sp(H)$-equivariant map.
The rational version
\[
\tau^{\Q} \colon \gr^{\Q}_\J \I \to \h^{\Q},
\]
where $\h^{\Q} = \h \otimes_\Z \Q$, is naturally an $\Sp(H^{\Q})$-map.
See \cite[(2.2.8)]{AsaNak}, \cite{Patzt} for instance.

\begin{rem}
The Lie algebra $\h$ of symplectic derivations was independently introduced by Kontsevich~\cite{Kon93} in the framework of formal symplectic geometry. 
\end{rem}

One could draw a parallel by saying that $\I$ is ``the group of IA-automorphisms that preserve $\zeta$'', and $\h$ is ``the Lie algebra of positive derivations that annihilate $\omega$.''   
We can think of Morita's refinement of $\tau$ as a ``linearization of the Dehn--Nielsen map.''
While $\rho$ is an isomorphism, neither $\tau$ nor $\tau^\Q$ is. 
In short, the Lie algebra $\h$ is still too large for the mapping class group, and this linearization process might be overlooking something important about self-homeomorphisms of the surface.
Therefore, it would be interesting to study topological interpretations for the cokernel of $\tau$.

We recall a diagrammatic description of the Lie algebra of symplectic derivations (see e.g., \cite{GaLe05}).
Since this is based on the aforementioned diagrammatic description of free Lie algebras, we work over the rationals.
Then, we have the following isomorphism of graded $\Q$-vector spaces:
\[
\h^\Q \overset{\cong}{\longleftarrow} \frac{\Q \{ \text{$H^{\Q}$-labeled planar trivalent trees} \}}{\text{(ML), (AS), (IHX)}}.
\]
Here, the grading on the right hand side is the number of trivalent vertices in the tree.
For example, the following figure represents an element of degree $3$.
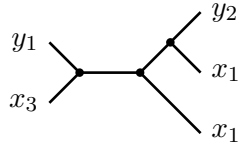
\[
\begin{tikzpicture}[x=4mm, y=4mm]
\draw (-2,0) -- (0,0) -- (1,1) -- (2,2) node[right]{$y_2$};
\draw (1,1) -- (2,0) node[right]{$x_1$};
\draw (0,0) -- (2,-2) node[right]{$x_1$};
\draw (-2,0) -- (-3,1) node[left]{$y_1$};
\draw (-2,0) -- (-3,-1) node[left]{$x_3$};
\fill (0,0) circle[radius=1.6pt];
\fill (1,1) circle[radius=1.6pt]; 
\fill (-2,0) circle[radius=1.6pt];
\end{tikzpicture}
\]
The above-mentioned isomorphism onto $\h^{\Q}$ is given as follows. 
First, note that the restriction to the degree-one part gives a graded $\Q$-linear isomorphism
\[
\Der^+(\L^\Q) \overset{\cong}{\to} \Hom(H^\Q, \L_{\ge 2}^{\Q}),
\quad D \mapsto D|_{H^{\Q}}.
\]
Here, $\L_{\ge 2}^{\Q} = \bigoplus_{k \ge 2} \L_k^{\Q}$.
Let $T$ be an $H^{\Q}$-labeled planar trivalent tree with $k$ trivalent vertices.
For each univalent vertex $v$ of $T$, let $x_v \in H^{\Q}$ be the label to $v$ and let $T_v$ be the $H^{\Q}$-labeled planar rooted trivalent tree obtained from $T$ by removing the label $x_v$, which we identify with an element in $\L_{k+1}$.
Then, $T$ corresponds to a symplectic derivation whose restriction to the degree-one part coincides with
\[
\bigg( y \mapsto \sum_{v} (x_v, y)\ T_v \bigg) \in \Hom(H^{\Q}, \L_{k+1}^{\Q}),
\]
where the sum is taken over all univalent vertices of $T$. 
Under this isomorphism, the Lie bracket of $\h^{\Q}$ is given by
\[
[T, T'] = \sum_{v,v'} (x_v,x_{v'}) 
\begin{tikzpicture}[baseline=-3pt, x=2mm, y=2mm]
\draw[ultra thick] (0,0) node[left=-3pt]{$T_v$} -- (3,0) node[right=-3pt]{$T'_{v'}$,};
\end{tikzpicture}
\]
where $(x_v,x_{v'})$ is the intersection number of the labels of $v$ and $v'$, and 
$
\begin{tikzpicture}[baseline=-3pt, x=2mm, y=2mm]
\draw[ultra thick] (0,0) node[left=-3pt]{$T_v$} -- (3,0) node[right=-3pt]{$T'_{v'}$};
\end{tikzpicture}
$
is the $H$-labeled planar trivalent tree obtained by gluing ``root-to-root'' $T_v$ and $T'_{v'}$.

\begin{example}
In degree one, we have $H^{\Q}$-labeled tripods on the diagrammatic side.
We have 
\[
\h^\Q(1) \overset{\cong}{\to} \wedge^3 H^{\Q}, 
\qquad
\begin{tikzpicture}[baseline=-2pt, x=4mm, y=4mm]
\fill (0,0) circle[radius=1.6pt];
\draw (0,0) -- (90:1.5) node[above]{$a$};
\draw (0,0) -- (210:1.5) node[left]{$b$};
\draw (0,0) -- (330:1.5) node[right]{$c$};
\end{tikzpicture}
\quad 
\mapsto 
\quad 
a \wedge b \wedge c.
\]
Note that this isomorphism holds over the integers: $\h(1) \cong \wedge^3 H$.
\end{example}

%%%%%%%
\subsection{Computations in degrees one and two}
We will make use of the diagrammatic expression of symplectic derivations by $H$-labeled planar trivalent trees.
(Note, however, that we no longer have an isomorphism between $\h$ and the corresponding $\Z$-module on the diagrammatic side.)

Johnson~\cite{J_aq} computed $\tau_1$.
Let $(\gamma, \delta)$ be a bounding pair which bounds a surface of genus one:
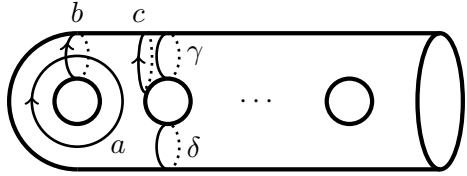
\[
\begin{tikzpicture}[x=3mm, y=3mm]
\draw[ultra thick] (0,0) circle[x radius=1, y radius=3];
\draw[ultra thick] (0,3) -- (-16,3);
\draw[ultra thick] (0,-3) -- (-16,-3);
\draw[ultra thick] (-16,3) arc[x radius=3, y radius=3, start angle=90, end angle=270];
\draw[ultra thick] (-4,0) circle[radius=1];
\draw (-8,0) node{$\cdots$};
\draw[ultra thick] (-12,0) circle[radius=1];
\draw[ultra thick] (-16,0) circle[radius=1];]
\draw (-12,3) arc[x radius=0.5, y radius=1, start angle=90, end angle=270];
\draw[dotted] (-12,1) arc[x radius=0.5, y radius=1, start angle=-90, end angle=90];
\draw (-12,-1) arc[x radius=0.5, y radius=1, start angle=90, end angle=270];
\draw[dotted] (-12,-3) arc[x radius=0.5, y radius=1, start angle=-90, end angle=90];
\draw (-12,2) node[right=3pt]{$\gamma$};
\draw (-12,-2) node[right=3pt]{$\delta$};
%%
%%% curve a %%%
\draw (-16,0) circle[radius=2];
\draw[->] (-18,0) -- (-18,0.2);
\draw (-14.2,-2) node{$a$};
%%% curve b %%%
\draw (-16,3) arc[x radius=0.5, y radius=1, start angle=90, end angle=270];
\draw[dotted] (-16,1) arc[x radius=0.5, y radius=1, start angle=-90, end angle=90];
\draw[->] (-16.45,2.4) -- (-16.45,2.6);
\draw (-16,3) node[above]{$b$};
%%% curve c %%%
\draw (-13,3) arc[x radius=0.3, y radius=1.3, start angle=90, end angle=270];
\draw[dotted] (-13,0.6) arc[x radius=0.3, y radius=1.3, start angle=-90, end angle=90];
\draw[->] (-13.3,1.8) -- (-13.3,2);
\draw (-13.3,3) node[above]{$c$};
%%%
\end{tikzpicture}
\]
Then, it holds that 
\[
\tau_1( t_{\gamma} t_{\delta}^{-1}) =
\begin{tikzpicture}[baseline=-2pt, x=4mm, y=4mm]
\fill (0,0) circle[radius=1.6pt];
\draw (0,0) -- (90:1.5) node[above]{$a$};
\draw (0,0) -- (210:1.5) node[left]{$b$};
\draw (0,0) -- (330:1.5) node[right]{$c$};
\end{tikzpicture}
\]
where $a$, $b$ and $c$ are the homology classes of oriented simple closed curves depicted in the figure.
Johnson proved that for $g\ge 2$ the map 
\[
\tau_1 \colon \gr_\J \I(1) \to \h(1) \cong \wedge^3 H
\]
is surjective and hence an isomorphism.
Furthermore, Johnson~\cite{J_Tor3} proved that $\tau_1\colon \I \to \wedge^3 H$ gives the abelianization of $\I$ modulo torsion for $g \ge 3$, namely it induces an isomorphism $\I^{\rm abel} \otimes_{\Z} \Q \overset{\cong}{\to} \wedge^3 H^\Q$.

The second Johnson homomorphism was computed by Morita~\cite[Section~1]{Mor89}.
Let $\gamma$ be a separating simple closed curve which bounds a subsurface of genus $h$:
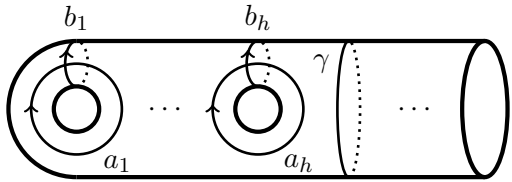
\[
\begin{tikzpicture}[x=3mm, y=3mm]
\draw[ultra thick] (0,0) circle[x radius=1, y radius=3];
\draw[ultra thick] (0,3) -- (-18,3);
\draw[ultra thick] (0,-3) -- (-18,-3);
\draw[ultra thick] (-18,3) arc[x radius=3, y radius=3, start angle=90, end angle=270];
\draw (-3,0) node{$\cdots$};
\draw[ultra thick] (-10,0) circle[radius=1];
\draw(-14,0) node{$\cdots$};
\draw[ultra thick] (-18,0) circle[radius=1];]
\draw (-6,3) arc[x radius=0.5, y radius=3, start angle=90, end angle=270];
\draw[dotted] (-6,-3) arc[x radius=0.5, y radius=3, start angle=-90, end angle=90];
\draw (-6,2) node[left=3pt]{$\gamma$};
%%
%%% curve a_1 %%%
\draw (-18,0) circle[radius=2];
\draw[->] (-20,0) -- (-20,0.2);
\draw (-16.2,-2.4) node{$a_1$};
%%% curve b_1 %%%
\draw (-18,3) arc[x radius=0.5, y radius=1, start angle=90, end angle=270];
\draw[dotted] (-18,1) arc[x radius=0.5, y radius=1, start angle=-90, end angle=90];
\draw[->] (-18.45,2.4) -- (-18.45,2.6);
\draw (-18,3) node[above]{$b_1$};
%%% curve a_h %%%
\draw (-10,0) circle[radius=2];
\draw[->] (-12,0) -- (-12,0.2);
\draw (-8.2,-2.4) node{$a_h$};
%%% curve b_h %%%
\draw (-10,3) arc[x radius=0.5, y radius=1, start angle=90, end angle=270];
\draw[dotted] (-10,1) arc[x radius=0.5, y radius=1, start angle=-90, end angle=90];
\draw[->] (-10.45,2.4) -- (-10.45,2.6);
\draw (-10,3) node[above]{$b_h$};
\end{tikzpicture}
\]
Morita proved the formula 
\[
\tau_2(t_{\gamma}) = \frac{1}{2} \sum_{i,j=1}^h 
\begin{tikzpicture}[baseline=-2pt, x=4mm, y=4mm]
\fill (0,0) circle[radius=1.6pt];
\fill (2,0) circle[radius=1.6pt];
\draw (0,0) -- (120:1.5) node[left]{$b_i$};
\draw (0,0) -- (240:1.5) node[left]{$a_i$};
\draw (0,0) -- (2,0);
\draw (2,0) --++(60:1.5) node[right]{$a_j$};
\draw (2,0) --++(300:1.5) node[right]{$b_j$};
\end{tikzpicture} 
\]
and showed that $\tau_2$ is an isomorphism modulo $2$-torsions and hence $\tau^{\Q}_2$ is an isomorphism:
\[
\tau^{\Q}_2 \colon \gr^{\Q}_\J \I(2) \overset{\cong}{\to} \h^{\Q}(2).
\]

\subsection{The Johnson cokernel problem for the mapping class group} \label{subsec:JCP}

As posed by Morita~\cite[Section~6, Main problem]{Mor93}, it is natural to formulate the Johnson cokernel problem for the mapping class group as follows.

\begin{problem}
Determine explicitly the cokernel of $\tau\colon \gr_{\J} \I \to \h$ and $\tau^{\Q}\colon \gr_{\J}^{\Q} \I \to \h^{\Q}$. 
\end{problem}

Let us confine ourselves to the rational version.
We have seen that $\tau^{\Q}_1$ and $\tau^{\Q}_2$ are surjective, so their cokernels are trivial.
However, in general, $\tau_k^{\Q}\colon \gr_{\J}^{\Q} \I(k) \to \h^{\Q}(k)$ is not surjective for $k \ge 3$.
Morita~\cite[Section~6]{Mor93} first found this fact.
He showed that for each odd $k \ge 3$ there is an $\Sp(H^{\Q})$-equivariant surjection $\h^{\Q}(k) \to S^k H^{\Q}$, called the trace map, which vanishes on the image of $\tau_k^{\Q}$.
Here, $S^k H^{\Q}$ is the $k$th symmetric power of $H^{\Q}$.
Morita's\index{Morita trace} trace\index{trace!Morita} map plays an important role in the determination of the cokernel of $\tau^{\Q}_3$, which was done by Asada and Nakamura~\cite{AsaNak} and Hain~\cite[Section~9]{Hain97}.
Furthermore, the cokernel of $\tau^{\Q}_4$ was determined by Morita~\cite[Proposition~6.1]{Mor99}.

As $k$ increases, it becomes more and more difficult to determine the cokernel of $\tau_k^{\Q}$.
The following result of Hain is very important.
The proof is deep; it uses the theory of mixed Hodge structures and the relative Malcev completion of the mapping class group.

\begin{thm}[Hain~\cite{Hain93, Hain97}, see also \cite{Hain_surv}] \label{thm:Hain}
If $g \ge 3$, the comparison map 
\[
c^{\Q} \colon \gr^{\Q}_{\LC} \I \to \gr^{\Q}_{\J} \I
\]
is surjective.
Hence the image of $\tau^{\Q}$ is generated by the degree one part $\h^{\Q}(1) \cong \wedge^3 H^{\Q}$. 
\end{thm}

With this theorem in hand, the Johnson cokernel problem for $\M$ can be, at least over the rationals, stated purely algebraically: it asks to determine the cokernel of the iterated Lie bracket $(\h^{\Q}(1))^{\otimes k} \to \h^{\Q}(k)$.
Morita, Sakasai and Suzuki have determined the explicit $\Sp(H^{\Q})$-irreducible decomposition of the cokernel of $\tau_k^{\Q}$ up to degree $9$ utilizing Hain's result, other techniques, and huge computer calculations. 
Their result can be found partially in \cite{MSS15}.

\begin{rem}
As we will mention in Section~\ref{subsec:Recent}, great advances have been made recently on the kernel of the comparison map $c^{\Q}$ and the Johnson cokernel problem for $\M$ in the stable range.
\end{rem}

One goal in the Johnson cokernel problem for $\M$ is to find a complete set of defining equations of the image of $\tau^{\Q}$ in the Lie algebra $\h^{\Q}$. 
We call such an equation an obstruction for the surjectivity of $\tau^{\Q}$; it is a nontrivial $\Sp(H^{\Q})$-map from $\h^{\Q}(k)$ (or a subspace of $\h^{\Q}(k)$) to an $\Sp(H^{\Q})$-module which vanishes on the image of $\tau_k^{\Q}$.
Morita's trace map is one example, and there are many others, some of which we will discuss later. 
Among others, there is a class of obstructions called the Galois obstructions, which come from the absolute Galois group ${\rm Gal}(\overline{\Q}/\Q)$ of the rational number field and which appear as trivial $\Sp(H^{\Q})$-representations.
See \cite{Nak96, Mat13, MSS26} for more details.

In the rest of this section, we explain an important obstruction for the surjectivity of $\tau^{\Q}$, which was introduced by Enomoto and Satoh~\cite{ES14}.
They considered the following $\Sp(H^{\Q})$-equivariant map: 
\[
\Tr^{\ES} \colon \h^{\Q}(k) \hookrightarrow \Der^k(\L) \overset{\Phi_k}{\to} \C_n(k).
\]
Namely, the Enomoto--Satoh\index{Enomoto--Satoh trace} trace\index{trace!Enomoto--Satoh} $\Tr^{\ES}$ is the restriction of the map $\Phi_k$, which is explained in Section~\ref{subsec:autFn}, to $\h^{\Q}(k)$.
Enomoto and Satoh showed that $\Tr^{\ES}$ vanishes on the image of $\tau^{\Q}_k$ and found many explicit components in the cokernel of $\tau^{\Q}$ which cannot be detected by Morita's trace map.
In terms of the diagrammatic presentation of $\h^{\Q}$, the map $\Tr^{\ES}$ is given as follows: 
\begin{align*}
\begin{tikzpicture}[baseline=0pt, x=4.5mm, y=4.5mm]
\fill (0,0) circle[radius=1.6pt];
\fill (3,0) circle[radius=1.6pt];
\fill (1.5,0) circle[radius=1.6pt];
\draw (0,0) -- (120:1) node[above=-2pt]{$a_2$};
\draw (0,0) -- (240:1) node[below=-2pt]{$a_1$};
\draw (0,0) -- (3,0);
\draw (1.5,0) -- (1.5,1) node[above=-2pt]{$b_3$};
\draw (3,0) --++(60:1) --++(0:1) node[right=-2pt]{$a_4$};
\draw (3,0) --++(60:1) --++(120:1) node[above=-2pt]{$b_3$};
\draw (3,0) --++(300:1) node[below=-2pt]{$b_1$};
\end{tikzpicture}
\quad 
& \overset{\Tr^{\ES}}{\mapsto}
\quad 
(a_1, b_1) \, 
\begin{tikzpicture}[baseline=0pt, x=4.5mm, y=4.5mm]
\fill (0,0) circle[radius=1.6pt];
\fill (3,0) circle[radius=1.6pt];
\fill (1.5,0) circle[radius=1.6pt];
\draw (0,0) -- (120:1) node[above=-2pt]{$a_2$};
\draw (0,0) -- (240:1);
\draw (0,0) -- (3,0);
\draw (1.5,0) -- (1.5,1) node[above=-2pt]{$b_3$};
\draw (3,0) --++(60:1) --++(0:1) node[right=-2pt]{$a_4$};
\draw (3,0) --++(60:1) --++(120:1) node[above=-2pt]{$b_3$};
\draw (3,0) --++(300:1);
\draw (0,0) -- (-0.5,-0.866025); 
\draw[dashed, ->-] (-0.5,-0.866025) to[bend right=60] (0,-1.5) -- (3,-1.5) to[bend right=60] (3.5,-0.866025); 
\end{tikzpicture}
\, + \cdots \\ 
&= \, |[b_3,a_4]b_3a_2| \pm |a_2 b_3[b_3,a_4]| + \cdots.
\end{align*}
To be a bit more precise, let $T$ be an $H^{\Q}$-labeled planar trivalent tree.
For two leaves of $T$, contract their labels by the intersection pairing of $H^{\Q}$ and add a dotted edge.
By using the IHX relation one can read from the result a cyclic word in two ways.
The sum of the two cyclic words thus obtained (with a suitable choice of sign) is the contribution from the two leaves we start with.
Up to a constant, $\Tr^{\ES}(T)$ is the sum of contributions from all pairs of two leaves of $T$.
See Conant~\cite[Theorem~4.2(2)]{Con15}. 

On the one hand, the Enomoto--Satoh trace is so powerful that it detects all the components in the cokernel of $\tau_k^{\Q}$ for $k \le 5$.
On the other hand, it is not sufficient to capture everything: it is known that there is a component in the cokernel of $\tau_6^{\Q}$ which cannot be detected by $\Tr^{\ES}$. 

There is a topological interpretation of the Enomoto--Satoh trace in terms of a certain loop operation on the surface, called the Turaev cobracket\index{Turaev cobracket} \cite{Tur91}.
It is defined by splitting free loops at their self-intersections: 
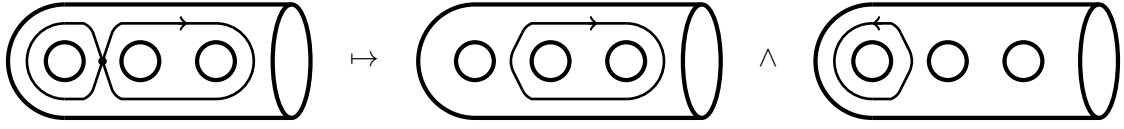
\[
\begin{tikzpicture}[baseline=-2pt, x=2.5mm, y=2.5mm]
\draw[ultra thick] (0,0) circle[x radius=1, y radius=3];
\draw[ultra thick] (0,3) -- (-12,3);
\draw[ultra thick] (0,-3) -- (-12,-3);
\draw[ultra thick] (-12,3) arc[x radius=3, y radius=3, start angle=90, end angle=270];
\draw[ultra thick] (-4,0) circle[radius=1];
\draw[ultra thick] (-8,0) circle[radius=1];
\draw[ultra thick] (-12,0) circle[radius=1];]
%%%
\draw (-12,2) arc[x radius=2, y radius=2, start angle=90, end angle=270];
\draw (-12,-2) -- (-11,-2) to[bend right=20] (-10.5,-1.5) -- (-9.5,1.5) to[bend left=20] (-9,2) -- (-4,2);
\draw[->] (-6,2) -- (-5.5,2);
\draw (-12,2) -- (-11,2) to[bend left=20] (-10.5,1.5) -- (-9.5,-1.5) to[bend right=20] (-9,-2) -- (-4,-2);
\draw (-4,-2) arc[x radius=2, y radius=2, start angle=-90, end angle=90];
%%%
\fill[ultra thick] (-10,0) circle[radius=1.6pt];
\end{tikzpicture}
\quad 
\mapsto 
\quad 
\begin{tikzpicture}[baseline=-2pt, x=2.5mm, y=2.5mm]
\draw[ultra thick] (0,0) circle[x radius=1, y radius=3];
\draw[ultra thick] (0,3) -- (-12,3);
\draw[ultra thick] (0,-3) -- (-12,-3);
\draw[ultra thick] (-12,3) arc[x radius=3, y radius=3, start angle=90, end angle=270];
\draw[ultra thick] (-4,0) circle[radius=1];
\draw[ultra thick] (-8,0) circle[radius=1];
\draw[ultra thick] (-12,0) circle[radius=1];]
%%%
\draw (-9.5,-1.5) to[bend right=20] (-9,-2) -- (-4,-2);
\draw (-4,-2) arc[x radius=2, y radius=2, start angle=-90, end angle=90];
\draw (-4,2) -- (-9,2) to[bend right=20] (-9.5,1.5);
\draw (-9.5,1.5) -- (-10,0.5) to[bend right=20] (-10,-0.5) -- (-9.5,-1.5);
%%%
\draw[->] (-6,2) -- (-5.5,2);
\end{tikzpicture}
\quad 
\wedge 
\quad
\begin{tikzpicture}[baseline=-2pt, x=2.5mm, y=2.5mm]
\draw[ultra thick] (0,0) circle[x radius=1, y radius=3];
\draw[ultra thick] (0,3) -- (-12,3);
\draw[ultra thick] (0,-3) -- (-12,-3);
\draw[ultra thick] (-12,3) arc[x radius=3, y radius=3, start angle=90, end angle=270];
\draw[ultra thick] (-4,0) circle[radius=1];
\draw[ultra thick] (-8,0) circle[radius=1];
\draw[ultra thick] (-12,0) circle[radius=1];]
%%%
\draw (-12,2) arc[x radius=2, y radius=2, start angle=90, end angle=270];
\draw (-12,2) -- (-11,2) to[bend left=20] (-10.5,1.5);
\draw (-12,-2) -- (-11,-2) to[bend right=20] (-10.5,-1.5);
\draw (-10.5,-1.5) -- (-10,-0.5) to[bend right=20] (-10,0.5) -- (-10.5,1.5);
\draw[->] (-11.5,2) -- (-12,2);
\end{tikzpicture}
\]
By construction, the Turaev cobracket is compatible with the action of the mapping class group on the $\Q$-span of the homotopy classes of free loops.
Using this fact, Kawazumi and Kuno~\cite[Corollary~6.3.3]{KK15} introduced an obstruction for the surjectivity of $\tau^{\Q}$, from which one can recover Morita's trace map \cite[Section~6.4]{KK15}. 
Furthermore, by introducing a framed version of the Turaev cobracket, Alekseev, Kawazumi, Kuno and Naef~\cite[Section~9]{AKKNg} refined this obstruction.
This refinement turns out to be equivalent to the Enomoto--Satoh trace. 
We only mention the key ingredients in their work: the Goldman--Turaev Lie bialgebra and its formality problem, the higher genus version of the Kashiwara--Vergne problem, and ideas from non-commutative Poisson geometry.

\subsection{$3$-dimensional topology}

There are well-known constructions in which the mapping class group interacts with $3$-dimensional topology: the mapping torus construction and Heegaard splittings.
We briefly summarize applications of the Johnson homomorphisms to invariants of $3$-manifolds.
In fact, the first abelian quotient of the Torelli group was given by Sullivan~\cite{Sull75} in the context of $3$-dimensional topology.
This can be seen as a precursor to the first Johnson homomorphism~\cite[p.\ 169]{J_surv}.

\subsubsection*{Massey products of mapping tori}
Given a mapping class $\varphi \in \M$, one can form the mapping torus\index{mapping torus} $T_{\varphi}$.
It is defined to be the $3$-manifold obtained from the Cartesian product $\Sigma \times [0,1]$ by identifying $(x,0)$ with $(\varphi(x),1)$.
Here we use the same letter $\varphi$ for a self-homeomorphism of $\Sigma$ which represents $\varphi$. 
The second projection endows $T_{\varphi}$ with the structure of a $\Sigma$-bundle over $S^1$.

Massey products\index{Massey product} are higher operations on the cohomology of a space. 
Johnson~\cite[p.\ 171]{J_surv} pointed out a relationship between the Johnson homomorphisms and Massey products of the mapping tori formed by elements of the Torelli group.
Kitano~\cite{Kit96} clarified this and proved that $\tau_k(\varphi)$ is equivalent to the Massey product of order $k$ of the mapping torus $T_{\varphi}$ for $\varphi \in \M(k)$.

\subsubsection*{Core of the Casson invariant}

For the moment, we denote $\Sigma = \Sigma_{g,1}$ to emphasize the genus and the number of boundary components of the surface.
We also use the symbols $\M_{g,1}$ and $\I_{g,1}$ for the corresponding mapping class group and Torelli group.
Let $\Sigma_g$ be a closed oriented surface of genus $g$.
We also fix an embedded closed disk $D \subset \Sigma_g$ and regard $\Sigma_{g,1} = \Sigma_g \setminus {\rm Int}\, D$.
We fix a standard embedding of $\Sigma_g$ into $\R^3$, and let $V_g$ be the genus $g$ handlebody\index{handlebody} defined as the bounded domain of $\R^3$ cut by $\Sigma_g$:
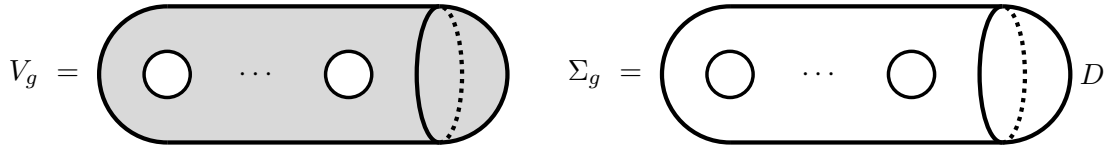
\[
V_g \ = \ 
\begin{tikzpicture}[baseline=-2pt, x=3mm, y=3mm]
%%%
\fill[gray!30] (-12,0) circle[radius=3];
\fill[gray!30] (-12,3) -- (-12,-3) -- (0,-3) -- (0,3) --cycle;
\fill[gray!30] (0,0) circle[radius=3];
%%%
\draw[ultra thick] (0,3) arc[x radius=1, y radius=3, start angle=90, end angle=270];
\draw[ultra thick, dotted] (0,-3) arc[x radius=1, y radius=3, start angle=-90, end angle=90];
\draw[ultra thick] (0,3) -- (-12,3);
\draw[ultra thick] (0,-3) -- (-12,-3);
\draw[ultra thick] (-12,3) arc[x radius=3, y radius=3, start angle=90, end angle=270];
\draw[ultra thick] (0,-3) arc[radius=3, start angle=-90, end angle=90];
\draw[ultra thick] (-4,0) circle[radius=1];
\draw (-8,0) node{$\cdots$};
\draw[ultra thick] (-12,0) circle[radius=1];
%%%%
\draw[fill=white] (-4,0) circle[radius=1];
\draw[fill=white] (-12,0) circle[radius=1];
\end{tikzpicture}
\qquad 
\Sigma_g \ = \ 
\begin{tikzpicture}[baseline=-2pt, x=3mm, y=3mm]
\draw[ultra thick] (0,3) arc[x radius=1, y radius=3, start angle=90, end angle=270];
\draw[ultra thick, dotted] (0,-3) arc[x radius=1, y radius=3, start angle=-90, end angle=90];
\draw[ultra thick] (0,3) -- (-12,3);
\draw[ultra thick] (0,-3) -- (-12,-3);
\draw[ultra thick] (-12,3) arc[x radius=3, y radius=3, start angle=90, end angle=270];
\draw[ultra thick] (0,-3) arc[radius=3, start angle=-90, end angle=90];
\draw[ultra thick] (-4,0) circle[radius=1];
\draw (-8,0) node{$\cdots$};
\draw[ultra thick] (-12,0) circle[radius=1];
%%%%
\draw[fill=white] (-4,0) circle[radius=1];
\draw[fill=white] (-12,0) circle[radius=1];
\draw (5,0) node[left]{$D$};
\end{tikzpicture}
\]
We can regard this as a Heegaard splitting of the $3$-sphere $S^3 = \R^3 \cup \{ \infty \}$: we have $S^3 = V_g \cup_{\iota} (-V_g)$, where $-V_g$ is a copy of $V_g$ with the opposite orientation and $\iota \colon \partial V_g = \Sigma_g \to \partial(-V_g) = \Sigma_g$ is a certain orientation preserving homeomorphism which restrict to the identity on $D$.

For a given element $\varphi \in \M_{g,1}$, one can form a $3$-manifold by twisting the aforementioned Heegaard splitting of $S^3$ by $\varphi$.
Namely, we define $M_\varphi := V_g \cup_{\iota \varphi} (-V_g)$, where $\varphi$ is extended to a self-homeomorphism of $\Sigma_g = \Sigma_{g,1} \cup D$ by using the identity map on $D$.
It is well known that any closed oriented $3$-manifold $M$ admits a decomposition of this form called a Heegaard splitting\index{Heegaard splitting} of $M$. 

A closed oriented $3$-manifold $M$ is called an (integral) homology $3$-sphere\index{homology $3$-sphere} if $H_*(M;\Z) \cong H_*(S^3;\Z)$.
If $\varphi \in \I_{g,1}$, then $M_\varphi$ is a homology $3$-sphere.
Conversely, any homology $3$-sphere can be obtained in this way for some $g$ and $\varphi \in \I_{g,1}$.
Let $\HS$ be the set of orientation-preserving homeomorphism classes of homology $3$-spheres.
It forms a commutative monoid under the connected sum operation whose identity element is the standard $3$-sphere $S^3$.
Thus we have a surjection
\[
\bigsqcup_{g} \I_{g,1} \to \HS, \quad \varphi \mapsto M_{\varphi}.
\]

Morita \cite[Proposition~2.3]{Mor89} proved that the pasting map $\varphi$ can be chosen from the Johnson kernel $\K = \M(2)$ to obtain any homology $3$-spheres, namely 
\[
\bigsqcup_{g} \M_{g,1}(2) \to \HS, \quad \varphi \mapsto M_{\varphi}
\]
is surjective.
Based on this result, Morita studied the Casson invariant\index{Casson invariant} $\lambda$, an important integer valued invariant for homology $3$-spheres which is of finite type in the sense of Ohtsuki~\cite{Oht96}, as a function on the Johnson kernel: $\lambda^* \colon \K \to \Z, \varphi \mapsto \lambda(M_{\varphi})$.
Morita showed that $\lambda^*$ is a homomorphism and decomposes into two homomorphisms.
One is determined by the second Johnson homomorphism $\tau_2$, and the other one, which Morita calls the core of the Casson invariant, is more mysterious and interesting.
This is a group homomorphism
\[
d \colon \K \to \Z
\]
which is invariant under the conjugation action of $\M$ and satisfies $d(t_{\gamma}) = 4h(h-1)$ if $t_{\gamma}$ is any BSCC-map of genus $h$.

Morita's homomorphism satisfies $d|_{\Gamma_3 \I} =0$ but $d|_{\M(3)} \neq 0$.
Hain~\cite{Hain93,Hain_surv} showed that for $g \ge 3$ there is an exact sequence
\[
0 \to \Q \to \gr_{\LC}^{\Q} \I(2) \overset{c^{\Q}}{\to} \gr_{\J}^{\Q} \I(2) \to 0.
\]
The homomorphism $d$ corresponds to the $1$-dimensional kernel of $c^{\Q}$.

\begin{rem}
Let $\K_g$ be the Johnson kernel for a closed surface of genus $g$, namely the subgroup of the mapping class group of $\Sigma_g$ generated by BSCC-maps.
Dimca, Hain and Papadima~\cite{DHP14} determined the rationalized abelianization of $\K_g$ for $g \ge 6$.
Their result was made more explicit by Morita, Sakasai and Suzuki~\cite{MSS20}, and Faes and Massuyeau~\cite{FM26} extended it for the Johnson kernel $\K$ of a once-bordered surface. 
Here, the core of the Casson invariant and a homomorphism $\tilde{\tau}_2 \colon \K \to H_3(\pi/\Gamma_3 \pi)$, which was introduced by Morita~\cite[Section~2]{Mor93} as a lift of the second Johnson homomorphism, play crucial roles.
\end{rem}

\begin{rem}
Morita's work on the Casson invariant has been developed by Garoufalidis and Levine~\cite{GaLe97, GaLe98} to the study of all finite type invariants of homology $3$-spheres in terms of the algebraic structure of the Torelli group.
\end{rem}

Finally, we mention results on the following natural question:
given $k$, is 
\[
\bigsqcup_{g} \M_{g,1}(k) \to \HS, \quad \varphi \mapsto M_{\varphi}
\]
surjective? 
The case $k = 2$ is affirmatively solved by Morita's result mentioned above.
The cases $k = 3$ and $k =4$ are affirmatively settled by Pitsch~\cite{Pit08} and Faes~\cite{Faes22}, respectively.
However, Pitsch and Riba~\cite{PitRiba} showed that this is no longer the case for $k = 5$.

%%%%
\subsubsection*{Homology cobordisms of surfaces}

There is a $3$-dimensional enlargement of the mapping class group introduced by Gousarrov~\cite{Gou99} and Habiro~\cite{Habiro00} and called the monoid $\mathcal{C}$ of homology cobordisms\index{homology cobordism} of $\Sigma$.
This monoid contains $\M$ as a submonoid, and it can be studied using common techniques and from common perspectives.
For example, the Johnson filtration $\mathcal{C} = \mathcal{C}(0) \supset \mathcal{C}(1) \supset \cdots$ and the Johnson homomorphism are naturally defined for $\mathcal{C}$.
On the Torelli part $\I\mathcal{C} = \mathcal{C}(1)$, there is another filtration called the $Y$-filtration, as a counterpart to the lower central series of $\I$.
A key feature is that the Johnson homomorphism for $\mathcal{C}$ is surjective, as was shown by Habegger~\cite{Hab00} and Garoufalidis and Levine~\cite{GaLe05}.
See also \cite{KM21} for an alternate proof. 

Since we limit ourselves to $\M$, we will not say much about $\C$.
For further details, we refer the reader to Habiro and Massuyeau's survey article~\cite{HabMas12} and references therein.

%%%%%
\subsection{Finite generation of the Johnson filtration}

As we have explained in Section~\ref{subsec:fg}, the Torelli group $\I = \M(1)$ for $g \ge 3$ is finitely generated.
However, it is not so when $g = 2$. 
This follows from the result of Mess~\cite{Mess} on the infinite generation of the Torelli group of a closed surface of genus $2$.

Mess's result also implies that for $g = 2$ the Johnson kernel $\K = \M(2)$ is not finitely generated. 
Whether $\K = \M(2)$ is finitely generated or not for $g \ge 3$ has received a lot of attention for some time.
See \cite[Problem~2.2(i)]{Mor99}.
A breakthrough was made by Ershov and He~\cite{EH18}, who proved that $\K = \M(2)$ is finitely generated for $g \ge 12$.
Furthermore, their result was soon refined with respect to the range of the genus and generalized to the higher terms:

\begin{thm}[Church--Ershov--Putman~\cite{CEP22}]
\begin{enumerate}
\item
Let $g \ge 4$.
Then, any subgroup of $\I$ containing $\Gamma_2 \I$ is finitely generated.
In particular, the Johnson kernel $\K$ is finitely generated.

\item 
For $k \ge 3$ and $g \ge 2k-1$, any subgroup of $\I$ containing $\Gamma_k \I$ is finitely generated.
In particular, $\M(k)$ is finitely generated.
\end{enumerate}
\end{thm}

%%%
\subsection{Torsion phenomena}

We have been working mostly over the rationals.
Here we collect classical and recent results on torsion phenomena around the Johnson filtration for the mapping class group.

Johnson~\cite{J_Tor3} determined the abelianization of the Torelli group $\I$.
He used the Birman--Craggs\index{Birman--Craggs homomorphism}  homomorphisms,\index{homomorphism!Birman--Craggs} which are $\Z_2$-valued homomorphisms from $\I$ related to the Rohlin invariant of homology $3$-spheres, to capture the torsion part of the abelianization $\I^{\rm abel} = \Gamma_1 \I/\Gamma_2 \I$.
As a natural extension, one can ask whether the higher successive quotients of the lower central series of the Torelli group have torsions or not. 
Faes, Massuyeau and Sato~\cite{FMS26} proved that $\Gamma_2 \I/\Gamma_3 \I$ is torsion free. 
Nozaki, Sato and Suzuki~\cite{NSS25} showed that for any odd $k$ and $g \ge 3k$, the quotient $\Gamma_k \I/\Gamma_{k+1} \I$ has nontrivial torsion elements.

As for the cokernel of the Johnson homomorphisms, the structure of the cokernel of $\tau_2$ (in particular its rank over $\Z_2$) was determined by Yokomizo~\cite{Yok02}.
See also Faes~\cite[Theorem~2.4]{Faes23} for a more conceptual description of Yokomizo's result.
Recently, Faes~\cite{Faes25} proves that for $k \ge 1$ and $g \ge k+2$ the cokernel of $\tau_{2k}$ has nontrivial $2$-torsion elements.

Nozaki, Sato and Suzuki~\cite{NSS22} proved that the Johnson kernel $\K_g$ for a closed surface of genus $g$ has nontrivial torsions.
This work is followed by a subsequent work of Faes and Massuyeau~\cite{FM26}, who give a lower bound for the cardinality of the torsions and prove the same for the Johnson kernel $\K = \K_{g,1}$ for a once-bordered surface.

%%%
\subsection{Recent progress on the Johnson cokernel problem} \label{subsec:Recent}

\subsubsection*{Stable injectivity of the comparison map}

As was mentioned in Theorem~\ref{thm:Hain}, the comparison map $c^{\Q} \colon \gr_{\LC}^{\Q} \I(k) \to \gr_{\J}^{\Q} \I(k)$ is surjective.
The injectivity of $c^{\Q}$ for $k \ge 3$ has been of interest (see e.g., \cite[Problem~6.2]{Mor99}).
The case $k=3$ is settled in \cite[Proposition~6.3]{Mor99}, and the cases $k = 4,5,6$ in \cite[Theorem~1.2 and Proposition~3.1]{MSS20}. 
Recently, there has been a significant progress on this question.

\begin{thm}[Kupers--Randal-Williams~\cite{KRW23} (modulo $\Sp(H^{\Q})$-invariants), Naef--Willwacher~\cite{NW26}] \label{thm:c_inj}
If $k \ge 3$ and $g \ge 3k$, the comparison map
\[
c^{\Q} \colon \gr_{\LC}^{\Q} \I(k) \to \gr_{\J}^{\Q} \I(k)
\]
is injective (and hence is an isomorphism).
\end{thm}

It was shown independently by Kupers and Randal-Williams~\cite{KRW23} and Felder, Naef and Willwacher~\cite{FNW23} that $\gr_{\LC}^{\Q} \I$ is stably Koszul (in the range $k \le g/3$).
As a consequence of this, Kupers and Randal-Williams~\cite[Theorem~B]{KRW23} proved that the kernel of $c^{\Q}$ consists of trivial $\Sp(H^{\Q})$-representations in the same stable range.
Naef and Willwacher~\cite[Theorem~1]{NW26} have made an improvement on this result by applications of their previous work \cite{FNW23} with Felder.
What is more, Naef and Willwacher~\cite[Theorem~2]{NW26} give a formula for the $\Sp(H^{\Q})$-irreducible decomposition of the cokernel of $\tau_k^{\Q}$ in a stable range.  
In the proof, a certain graph complex is used to reformulate the Johnson homomorphism.
Notably, they introduce a series of maps called the higher Enomoto--Satoh traces, defined on a series of descending Lie subalgebras of (a graphical model of) $\h^{\Q}$. 
 
\subsubsection*{Conant's stable description}

Let $\cok(k)$ be the cokernel of the $k$th Johnson homomorphism $\tau_k^{\Q}$.
Theorem~\ref{thm:Hain} says that $\cok(k)$ coincides with the cokernel of the iterated Lie bracket $(\h^{\Q}(1))^{\otimes k} \to \h^{\Q}(k)$.
Based on this fact and using the hairy graph complex~\cite{CKV13, CKV15}, Conant~\cite{Con16} introduced a graphical model which describes $\cok(k)$ in a stable range $g \gg k$.
Let us give a few more details.

By generalities on the representation theory of the symplectic group we have 
\[
\cok(k) = \cok_{k,\{ k\}} \oplus \cok_{k,\{k-2\}} \oplus \cdots
\]
where $\cok_{k,\{l\}}$ is a direct sum of irreducible $\Sp(H^{\Q})$-submodules of $\cok(k)$ corresponding to Young diagrams of size $l$.
Conant defined a bigraded space $\widetilde{\Omega} = \bigoplus_{r \ge 1} \widetilde{\Omega}_r = \bigoplus_{r \ge 1} \bigoplus_{m \ge 0} \widetilde{\Omega}_{r,m}$, where $\widetilde{\Omega}_{r,m}$ is $\Q$-linearly spanned by planar trivalent trees whose $2r$ leaves are paired by $r$ dotted edges and the remaining $m$ leaves are labeled by elements of $H^{\Q}$ (thus there are $2r+m$ leaves).
Such trees are called hairy Lie graphs.
The following figure is an example for $(r,m) = (2,6)$.
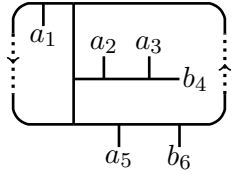
\[
\begin{tikzpicture}[baseline=-3pt, x=2mm, y=2mm]
%%%
%%% A hairy Lie graph %%%
%%%
%%% middle vertical line %%%
\draw (0,4) -- (0,-4);  
%%% right dotted edge %%%
\draw[very thick, dotted, ->-] (10,-2) -- (10,2);   
%%% left dotted edge %%%
\draw[very thick, dotted, ->-] (-4,2) -- (-4,-2); 
%%% big boundaries %%%
\draw (0,4) to (-3,4) to[bend right=30] (-4,3) to (-4,2);
\draw (0,-4) to (-3,-4) to[bend left=30] (-4,-3) to (-4,-2);
\draw (0,4) to (9,4) to[bend left=30] (10,3) to (10,2);
\draw (0,-4) to (9,-4) to[bend right=30] (10,-3) to (10,-2);
%%% hairs %%%
%% upper left %%
\draw (-2,4) -- (-2,2.5) node[below=-3pt]{$a_1$};
%% middle horizontal %%
\draw (0,-1) -- (7,-1) node[right=-3pt]{$b_4$};
\draw (2,-1) -- (2,0.5) node[above=-3pt]{$a_2$};
\draw (5,-1) -- (5,0.5) node[above=-3pt]{$a_3$};
%% lower right %%
\draw (3,-4) -- (3,-5.5) node[below=-3pt]{$a_5$};
\draw (7,-4) -- (7,-5.5) node[below=-5pt]{$b_6$};
\end{tikzpicture}
\]
In $\widetilde{\Omega}$, we have the relations (ML), (AS) and (IHX) appearing in the diagrammatic description of $\h^{\Q}$, the relation that flipping direction of a dotted edge gives a minus sign, and finally a certain relation corresponding to the image of the iterated Lie bracket $(\h^{\Q}(1))^{\otimes k} \to \h^{\Q}(k)$.
Then the trace map by Conant, Kassabov and Vogtmann~\cite{CKV13} induces a map
\[
\bigoplus_{r = 1}^{\lfloor k/2 \rfloor} \widetilde{\Tr}_r \colon
\h^{\Q}(k) \to \bigoplus_{r = 1}^{\lfloor k/2 \rfloor} \widetilde{\Omega}_{r,k+2-2r}.
\] 
For an $H^{\Q}$-labeled planar trivalent tree $T$, the  $r$-loop trace $\widetilde{\Tr}_r$ adds $r$ dotted edges to $T$ using contraction by the intersection pairing in all possible unordered ways.
It follows by construction that $\widetilde{\Tr}_r$ vanishes on the image of $\tau_k^{\Q}$, and the $1$-loop trace $\widetilde{\Tr}_1$ is essentially the same as the Enomoto--Satoh trace. 

Let $(H^{\Q})^{\langle m \rangle}$ be the intersection of the kernels of the pairwise contractions $(H^{\Q})^{ \otimes m} \to (H^{\Q})^{ \otimes (m-2)}$ by the intersection pairing.
Let $\widetilde{\Omega}_{r,\langle m \rangle}$ be the subspace of $\widetilde{\Omega}_{r,m}$ spanned by hairy Lie graphs whose leaves are labeled by tensors in $H^{\langle m \rangle}$.
Composing $\widetilde{\Tr}_r$ with natural projections $\widetilde{\Omega}_{r,m} \to \widetilde{\Omega}_{r,\langle m \rangle}$ we obtain
\[
\bigoplus_{r = 1}^{\lfloor k/2 \rfloor} \widetilde{\Tr}_r \colon
\h^{\Q}(k)/{\rm Im}\, \tau_k^{\Q} \to \bigoplus_{r = 1}^{\lfloor k/2 \rfloor} \widetilde{\Omega}_{r,\langle k+2-2r \rangle},
\] 
and for $g\gg k$ it is an $\Sp(H)$-isomorphism so that $\widetilde{\Omega}_{r, \langle k+2-2r \rangle} = \cok_{k,\{ k+2-2r\}}$.

It seems difficult to work directly with the last defining relation for $\widetilde{\Omega}_r$ for general $r$.
Recently, the case $r=2$ has been investigated by Kuno and Sato~\cite{KS25}.
They obtained an explicit presentation of $\widetilde{\Omega}_2$, and used it to capture certain components in the Johnson cokernel which cannot be detected by the Enomoto--Satoh trace.

It would be interesting to look for topological meanings for the higher loop traces $\widetilde{\Tr}_r$ for $r \ge 2$, and their relationships with the higher Enomoto--Satoh traces introduced by Naef and Willwacher~\cite{NW26}.


\begin{thebibliography}{99}

\bibitem{AKKNg}A.\ Alekseev, N.\ Kawazumi, Y.\ Kuno and F.\ Naef, 
The Goldman--Turaev Lie bialgebra and the Kashiwara--Vergne problem in higher genera,
arXiv:1804.09566v3, to appear in Mem.\ Amer.\ Math.\ Soc.

\bibitem{Andre}S.\ Andreadakis, 
On the automorphisms of free groups and free nilpotent groups, 
Proc.\ London Math.\ Soc. (3) {\bf 15} (1965), 239--268.

\bibitem{AsaNak}M.\ Asada and H.\ Nakamura, 
On graded quotient modules of mapping class groups of surfaces,
Israel J.\ Math. {\bf 90} (1995), no.~1--3, 93--113.

\bibitem{Bach}S.\ Bachmuth,
Induced automorphisms of free groups and free metabelian groups,
Trans.\ Amer.\ Math.\ Soc. {\bf 122} (1966), 1--17.

\bibitem{Barth}L.\ Bartholdi, 
Automorphisms of free groups I,
New York J.\ Math. {\bf 19} (2013), 395--421.

\bibitem{BarthEr}L.\ Bartholdi, 
Automorphisms of free groups I--erratum,
New York J.\ Math. {\bf 22} (2016), 1135--1137.

\bibitem{CEP22}T.\ Church, M.\ Ershov and A.\ Putman, 
On finite generation of the Johnson filtrations,
J.\ Eur.\ Math.\ Soc. (JEMS) {\bf 24} (2022), no.~8, 2875--2914.

\bibitem{CHP11}F.\ Cohen, A.\ Heap and A.\ Pettet, 
On the Andreadakis--Johnson filtration of the automorphism group of a free group,
J.\ Algebra {\bf 329} (2011), 72--91.

\bibitem{CohPak}F.\ Cohen and J.\ Pakianathan,
On automorphism groups of free groups, and their nilpotent quotients,
preprint.

\bibitem{Con15}J.\ Conant, 
The Johnson cokernel and the Enomoto--Satoh invariant, 
Algebr.\ Geom.\ Topol. {\bf 15} (2015), no.~2, 801--821.

\bibitem{Con16}J.\ Conant, 
Addendum to ``The Johnson cokernel and the Enomoto--Satoh invariant'': the
ES-trace detects all top-level partitions, arXiv:1610.05220.

\bibitem{CKV13}J.\ Conant, M.\ Kassasbov and K.\ Vogtmann, 
Hairy graphs and the unstable homology of ${\rm Mod}(g,s)$, ${\rm Out}(F_n)$ and ${\rm Aut}(F_n)$, 
J.\ Topol. {\bf 6} (2013), no.~1, 119--153.

\bibitem{CKV15}J.\ Conant, M.\ Kassasbov and K.\ Vogtmann, 
Higher hairy graph homology,
Geom.\ Dedicata {\bf 176} (2015), 345--374.

\bibitem{Dar19}J.\ Darn\'{e},
On the stable Andreadakis problem,
J.\ Pure Appl.\ Algebra {\bf 223} (2019), no.~12, 5484--5525.

\bibitem{Dar21}J.\ Darn\'{e}, 
On the Andreadakis problem for subgroups of $\IA_n$,
Int.\ Math.\ Res.\ Not.\ IMRN {\bf 2021}, no.~19, 14720--14742.

\bibitem{Dar25}J.\ Darn\'{e}, 
Braids, inner automorphisms and the Andreadakis problem,
Ann.\ Inst.\ Fourier (Grenoble) {\bf 75} (2025), no.~1, 185--227.

\bibitem{DES}J.\ Darn\'{e}, N.\ Enomoto and T.\ Satoh, 
The Andreadakis Problem for the McCool groups, arXiv:2502.20652.

\bibitem{DHP14}A.\ Dimca, R.\ Hain and S.\ Papadima, 
The abelianization of the Johnson kernel, 
J.\ Eur.\ Math.\ Soc. (JEMS) {\bf 16} (2014), no.~4, 805--822.

\bibitem{Dri91}V.\ G.\ Drinfeld, 
On quasitriangular quasi-Hopf algebras and on a group that is closely connected with ${\rm Gal}(\overline{\Q}/\Q)$,
Leningrad Math.\ J. {\bf 2} (1991) no.~4, 829--860.

\bibitem{ES14}N.\ Enomoto and T.\ Satoh, 
New series in the Johnson cokernels of the mapping class groups of surfaces,
Algebr.\ Geom.\ Topol. {\bf 14} (2014), no.~2, 627--669.

\bibitem{EH18}M.\ Ershov and S.\ He, 
On finiteness properties of the Johnson filtrations, 
Duke Math.\ J. {\bf 167} (2018), no.~9, 1713--1759.

\bibitem{Faes22}Q.\ Faes,
Triviality of the $J_4$-equivalence among homology 3-spheres,
Trans.\ Amer.\ Math.\ Soc. {\bf 375} (2022), no.~9, 6597--6620.

\bibitem{Faes23}Q.\ Faes, 
The handlebody group and the images of the second Johnson homomorphism, 
Algebr.\ Geom.\ Topol. {\bf 23} (2023), no.~1, 243--293.

\bibitem{Faes25}Q.\ Faes, 
Torsion in the cokernels of the Johnson homomorphisms,
Int.\ Math.\ Res.\ Not.\ IMRN {\bf 2025}, no.~10, Paper No. rnaf121, 32 pp.

\bibitem{FM26}Q.\ Faes and G.\ Massuyeau, 
On the non-triviality of the torsion subgroup of the abelianized Johnson kernel,
Ann.\ Inst.\ Fourier (Grenoble) {\bf 76} (2026), no.~1, 425--475.

\bibitem{FMS26}Q.\ Faes, G.\ Massuyeau and M.\ Sato, 
On the degree-two part of the associated graded of the lower central series of the Torelli group,
Math.\ Proc.\ Cambridge Philos.\ Soc. {\bf 180} (2026), no.~2, 303--341.

\bibitem{Farb}B.\ Farb,
Automorphisms of $F_n$ which act trivially on homology, preprint.

\bibitem{FM12}B.\ Farb and D.\ Margalit, 
A primer on mapping class groups, 
Princeton Math.\ Ser., {\bf 49}, 
Princeton University Press, Princeton, NJ, 2012. xiv+472 pp.

\bibitem{FNW23}M.\ Felder, F.\ Naef and T.\ Willwacher, 
Stable cohomology of graph complexes,
Selecta Math.\ (N.S.) {\bf 29} (2023), no.~2, Paper No.~23, 72 pp.

\bibitem{For90}E.\ Formanek, 
Characterizing a free group in its automorphism group, 
J.\ Algebra {\bf 133} (1990), no.~2, 424--432.

\bibitem{GaLe97}S.\ Garoufalidis and J.\ Levine, 
Finite type $3$-manifold invariants, the mapping class group and blinks,
J.\ Differential Geom. {\bf 47} (1997), no.~2, 257--320.

\bibitem{GaLe98}S.\ Garoufalidis and J.\ Levine, 
Finite type $3$-manifold invariants and the structure of the Torelli group I,
Invent.\ Math. {\bf 131} (1998), no.~3, 541--594.

\bibitem{GaLe05}S.\ Garoufalidis and J.\ Levine, 
Tree-level invariants of three-manifolds, Massey products and the Johnson homomorphism,
Graphs and patterns in mathematics and theoretical physics, 173--203.
Proc.\ Sympos.\ Pure Math., {\bf 73},
American Mathematical Society, Providence, RI, 2005.

\bibitem{Gou99}M.\ Goussarov,
Finite type invariants and $n$-equivalence of $3$-manifolds,
C.\ R.\ Acad.\ Sci.\ Paris S\'{e}r.\ I Math. {\bf 329} (1999), no.~6, 517--522.

\bibitem{Hab00}N.\ Habegger, 
Milnor, Johnson, and tree level perturbative invariants, preprint (2000).

\bibitem{Habiro00}K.\ Habiro, 
Claspers and finite type invariants of links,
Geom.\ Topol. {\bf 4} (2000), no.~1, 1--83.

\bibitem{HabMas12}K.\ Habiro and G.\ Massuyeau, 
From mapping class groups to monoids of homology cobordisms: a survey,
Handbook of Teichm\"{u}ller theory, Vol.~III, 465--529.
IRMA Lect.\ Math.\ Theor.\ Phys., {\bf 17}, 
European Mathematical Society (EMS), Z\"{u}rich, 2012.

\bibitem{Hain93}R.\ Hain, 
Completions of mapping class groups and the cycle $C - C^{-}$,
Mapping class groups and moduli spaces of Riemann surfaces (G\"{o}ttingen, 1991/Seattle, WA, 1991), 75--105, 
Contemp.\ Math. {\bf 150}
American Mathematical Society, Providence, RI, 1993.

\bibitem{Hain97}R.\ Hain, 
Infinitesimal presentations of the Torelli groups, 
J.\ Amer.\ Math.\ Soc. {\bf 10} (1997), no.~3, 597--651.

\bibitem{Hain15}R.\ Hain, 
Genus 3 mapping class groups are not K\"{a}hler, 
J.\ Topol. {\bf 8} (2015), no.~1, 213--246.

\bibitem{Hain_surv}R.\ Hain, 
Johnson homomorphisms,
EMS Surv.\ Math.\ Sci. {\bf 7} (2020), no.~1, 33--116.

\bibitem{J_aq}D.\ Johnson, 
An abelian quotient of the mapping class group $\I_g$, 
Math.\ Ann. {\bf 249} (1980), no.~3, 225--242.

\bibitem{J_surv}D.\ Johnson, 
A survey of the Torelli group,
Low-dimensional topology (San Francisco, Calif., 1981), 165--179.
Contemp.\ Math., {\bf 20},
American Mathematical Society, Providence, RI, 1983.

\bibitem{J_Tor1}D.\ Johnson, 
The structure of the Torelli group I: A finite set of generators for $\mathcal{I}$, 
Ann.\ of Math. (2) {\bf 118} (1983), no.~3, 423--442.

\bibitem{J_Tor2}D.\ Johnson, 
The structure of the Torelli group II: A characterization of the group generated by twists on bounding curves,
Topology {\bf 24} (1985), no.~2, 113--126.

\bibitem{J_Tor3}D.\ Johnson, 
The structure of the Torelli group III: The abelianization of $\mathcal{I}$, 
Topology {\bf 24} (1985), no.~2, 127--144.

\bibitem{Kaw05}N.\ Kawazumi,
Cohomological aspects of Magnus expansions, arXiv:math/0505497.

\bibitem{KK15}N.\ Kawazumi and Y.\ Kuno, 
Intersection of curves on surfaces and their applications to mapping class groups, 
Ann.\ Inst.\ Fourier (Grenoble) {\bf 65} (2015), no.~6, 2711--2762.

\bibitem{KK16}N.\ Kawazumi and Y.\ Kuno, 
The Goldman--Turaev Lie bialgebra and the Johnson homomorphisms,
Handbook of Teichm\"{u}ller theory, Vol.~V, 97--165.
IRMA Lect.\ Math.\ Theor.\ Phys., {\bf 26}, 
European Mathematical Society (EMS), Z\"{u}rich, 2016.

\bibitem{Kit96}T.\ Kitano, 
Johnson's homomorphisms of subgroups of the mapping class group, the Magnus expansion and Massey higher products of mapping tori,
Topology Appl. {\bf 69} (1996), no.~2, 165--172.

\bibitem{Koh85}T.\ Kohno, 
S\'{e}rie de Poincar\'{e}--Koszul associ\'{e}e aux groupes de tresses pures,
Invent.\ Math. {\bf 82} (1985), no.~1, 57--75.

\bibitem{Kon93}M.\ Kontsevich, 
Formal (non)commutative symplectic geometry, The Gel'fand Mathematical Seminars, 1990--1992, 173--187.
Birkh\"{a}user Boston, Inc., Boston, MA, 1993.

\bibitem{KM21}Y.\ Kuno and G.\ Massuyeau, 
Generalized Dehn twists on surfaces and homology cylinders, 
Algebr.\ Geom.\ Topol. {\bf 21} (2021), no.~2, 697--754.

\bibitem{KS25}Y.\ Kuno and M.\ Sato, 
On the 2-loop part of the Johnson cokernel,
arXiv:2508.19041v2.

\bibitem{KRW23}A.\ Kupers and O.\ Randal-Williams, 
On the Torelli Lie algebra,
Forum Math.\ Pi {\bf 11} (2023), Paper No.~e13, 47 pp.

\bibitem{Lab70}J.\ Labute, 
On the descending central series of groups with a single defining relation,
J.\ Algebra {\bf 14} (1970), 16--23.

\bibitem{Laz54}M.\ Lazard,
Sur les groupes nilpotents et les anneaux de Lie,
Ann.\ Sci.\ \'{E}cole Norm.\ Sup. (3) {\bf 71} (1954), no.~2, 101--190.

\bibitem{Levine01Add}J.\ Levine, 
Addendum and correction to: ``Homology cylinders: an enlargement of the mapping class group,''
Algebr.\ Geom.\ Topol. {\bf 2} (2002), 1197--1204.

\bibitem{MKS}W.\ Magnus, A.\ Karrass and D.\ Solitar, 
Combinatorial group theory,
Presentations of groups in terms of generators and relations.
Reprint of the 1976 second edition.
Dover Publications, Inc., Mineola, NY, 2004. xii+444 pp.

\bibitem{MasSak}G.\ Massuyeau and T.\ Sakasai,
Morita's trace maps on the group of homology cobordisms,
J.\ Topol.\ Anal. {\bf 12} (2020), no.~3, 775--818.

\bibitem{Mat13}M.\ Matsumoto, 
Introduction to arithmetic mapping class groups,
Moduli spaces of Riemann surfaces, 319--356.
IAS/Park City Math.\ Ser., {\bf 20},
American Mathematical Society, Providence, RI, 2013.

\bibitem{Mess}G.\ Mess,
The Torelli groups for genus 2 and 3 surfaces,
Topology {\bf 31} (1992), no.~4, 775--790.

\bibitem{Mor89}S.\ Morita, 
Casson's invariant for homology 3-spheres and characteristic classes of surface bundles I, 
Topology {\bf 28} (1989), no.~3, 305--323.

\bibitem{Mor91}S.\ Morita, 
On the structure of the Torelli group and the Casson invariant,
Topology {\bf 30} (1991), no.~4, 603--621.

\bibitem{Mor93}S.\ Morita, 
Abelian quotients of subgroups of the mapping class group of surfaces, 
Duke Math.\ J., {\bf 70} (1993), 699--726.

\bibitem{Mor99}S.\ Morita, 
Structure of the mapping class groups of surfaces: a survey and a prospect,
Proceedings of the Kirbyfest (Berkeley, CA, 1998), 349--406.
Geom.\ Topol.\ Monogr., {\bf 2}, 
Geometry \& Topology Publications, Coventry, 1999.

\bibitem{MSS15}S.\ Morita, T.\ Sakasai and M.\ Suzuki, 
Structure of symplectic invariant Lie subalgebras of symplectic derivation Lie algebras,
Adv.\ Math. {\bf 282} (2015), 291--334.

\bibitem{MSS20}S.\ Morita, T.\ Sakasai and M.\ Suzuki, 
Torelli group, Johnson kernel, and invariants of homology spheres,
Quantum Topol. {\bf 11} (2020), no.~2, 379--410.

\bibitem{MSS26}S.\ Morita, T.\ Sakasai and M.\ Suzuki, 
The first Galois obstruction in the Johnson cokernel,
arXiv:2608.17673.

\bibitem{NW26}F.\ Naef and T.\ Willwacher, 
The Johnson homomorphism, embedding calculus and graph complexes,
arXiv:2602.09915v2.

\bibitem{Nak96}H.\ Nakamura, 
Coupling of universal monodromy representations of Galois--Teichm\"{u}ller modular groups, 
Math.\ Ann.\ {\bf 304} (1996), no.~1, 99--119.

\bibitem{NSS22}Y.\ Nozaki, M.\ Sato and M.\ Suzuki, 
Abelian quotients of the $Y$-filtration on the homology cylinders via the LMO functor,
Geom.\ Topol. {\bf 26} (2022), no.~1, 221--282.

\bibitem{NSS25}Y.\ Nozaki, M.\ Sato and M.\ Suzuki, 
Torsion elements in the associated graded of the $Y$-filtration of the monoid of homology cylinders,
J.\ Topol. {\bf 18} (2025), no.~3, Paper No.~e70028, 30 pp.

\bibitem{Oht96}T.\ Ohtsuki, 
Finite type invariants of integral homology 3-spheres,
J.\ Knot Theory Ramifications {\bf 5} (1996), no.~1, 101--115.

\bibitem{Pa01}S.\ Papadima, 
Braid commutators and homogenous Campbell--Hausdorff tests,
Pacific J.\ Math. {\bf 197} (2001), no.~2, 383--416. 

\bibitem{PaSu12}S.\ Papadima and A.\ Suciu, 
Homological finiteness in the Johnson filtration of the automorphism group of a free group, 
J.\ Topol. {\bf 5} (2012), no.~4, 909--944. 

\bibitem{Patzt}P.\ Patzt, 
Representation stability for filtrations of Torelli groups, 
Math.\ Ann.\ {\bf 372} (2018), no.~1, 257--298.

\bibitem{Pettet}A.\ Pettet,
The Johnson homomorphism and the second cohomology of $\IA_n$,
Algebr.\ Geom.\ Topol. {\bf 5} (2005), no.~2, 725--740.

\bibitem{Pit08}W.\ Pitsch, 
Integral homology 3-spheres and the Johnson filtration,
Trans.\ Amer.\ Math.\ Soc. {\bf 360} (2008), no.~6, 2825--2847.

\bibitem{PitRiba}W.\ Pitsch and R.\ Riba, 
Finite type invariants in low degrees and the Johnson filtration, 
arXiv:2311.09924v4.

\bibitem{Put12}A.\ Putman, 
Small generating sets for the Torelli group, 
Geom.\ Topol. {\bf 16} (2012), no.~1, 111--125.

\bibitem{Put18}A.\ Putman, 
The Johnson homomorphism and its kernel, 
J.\ Reine Angew.\ Math. {\bf 735} (2018), 109--141.

\bibitem{Satoh06}T.\ Satoh, 
New obstructions for the surjectivity of the Johnson homomorphism of the automorphism group of a free group,
J.\ London Math.\ Soc.\ (2) {\bf 74} (2006), no.~2, 341--360.

\bibitem{Satoh12}T.\ Satoh, 
On the lower central series of the IA-automorphism group of a free group,
J.\ Pure Appl.\ Algebra {\bf 216} (2012), no.~3, 709--717.

\bibitem{Satoh12MC}T.\ Satoh, 
On the Johnson filtration of the basis-conjugating automorphism group of a free group,
Michigan Math.\ J. {\bf 61} (2012), no.~1, 87--105.

\bibitem{Satoh16}T.\ Satoh, 
A survey of the Johnson homomorphisms of the automorphism groups of free groups and related topics,
Handbook of Teichm\"{u}ller theory, Vol.\ V, 167--209.
IRMA Lect.\ Math.\ Theor.\ Phys.\ {\bf 26},
European Mathematical Society (EMS), Z\"{u}rich, 2016.

\bibitem{Satoh17}T.\ Satoh,
On the Andreadakis conjecture restricted to the ``lower-triangular'' automorphism groups of free groups,
J.\ Algebra Appl. {\bf 16} (2017), no.~5, 1750099, 31 pp.

\bibitem{Satoh19}T.\ Satoh,
The third subgroup of the Andreadakis--Johnson filtration of the automorphism group of a free group,
J.\ Group Theory {\bf 22} (2019), no.~1, 41--61.

\bibitem{Sull75}D.\ Sullivan, 
On the intersection ring of compact three manifolds,
Topology 14 (1975), no.~3, 275--277.

\bibitem{Tur91}V.\ G.\ Turaev, 
Skein quantization of Poisson algebras of loops on surfaces, 
Ann.\ Sci.\ \'{E}cole Norm.\ Sup. (4) {\bf 24} (1991), no.~6, 635--704.

\bibitem{Yok02}Y.\ Yokomizo,  
An ${\rm Sp}(2g;\Z_2)$-module structure of the cokernel of the second Johnson homomorphism,
Topology Appl. {\bf 120} (2002), no.~3, 385--396.

\end{thebibliography}
\end{document}